\documentclass{amsart}

\usepackage[all]{xy} \usepackage{amssymb} \usepackage{amsmath} \usepackage{graphicx} \usepackage[english]{babel} \usepackage{color}
\usepackage{epsfig}
\usepackage{appendix}

\theoremstyle{definition} 
\newtheorem{definition}{Definition} 
\newtheorem{remark}[definition]{Remark}
\newtheorem{example}[definition]{Example} 
\newtheorem{acknowledgement}{Acknowledegment} 
\newtheorem{notation}[definition]{Notation}
\theoremstyle{plain} 
\newtheorem{lemma}[definition]{Lemma} 
\newtheorem{proposition}[definition]{Proposition}
\newtheorem{corollary}[definition]{Corollary}
 \newtheorem{theorem}[definition]{Theorem}

\def\SS{\mathrm{\sqcup\!\sqcup^\star\!~}} 
\def\Ss{\mathrm{\sqcup\!\sqcup\!~}}
 \def\DD{\Delta_\prec}
  \def\Diag{\mathfrak D}
  \def\Diagto{\stackrel{\to}{\Diag}} 
\def\oo{\omega} 
\def\mn{{}^{-1}} \def\t{\otimes}

\DeclareMathOperator{\sgn}{sgn}

\newcommand{\Vtm}[1]{V_m^{\otimes #1}} 
\newcommand{\Vt}[1]{V^{\otimes #1}} 

\newcommand{\G}[1]{{}^{\cdot #1}\!} 
 \newcommand{\GS}[2]{{}^{#1 #2}\!\!}

 \def\NN{{\mathbb{N}}}

\def\KK{{\mathbb{K}}}

\def\Sp{\mathit{sp}} 
\def\PP{\mathcal{P}}

\def\GG{\mathcal{G}}
\def\GGto{\stackrel{\to}{\GG}} 
\def\GGG{\mathfrak {G}} 
\def\HH{\mathcal{H}} 
\def\gg{\mathfrak{g}} 
\def\hh{\mathfrak{h}}

\def\dim{\mathrm{dim\, }} 
\def\Prim{\mathrm{Prim\, }} 
\def\Id{\mathrm{Id }} 

 \def\Com{\mathrm{Com}}
\def\std{\mathrm{std}}
 \def\min{\mathrm{min}} 
 \def\max{\mathrm{max}} 
 \def\norm{\mathrm{norm}}
  \def\Im{\mathop{\rm Im}}

 \def\PBT{\mathrm{PBT}}

 \newcommand{\commentc}[1]{} 

\begin{document}

\author[E.Burgunder]{Emily Burgunder} \address{ Institut de Math\'ematiques et de mod\'elisation de Montpellier \\ UMR CNRS 5149\\
D\'epartement de math\'ematiques\\ Universit\'e Montpellier II\\ Place Eug\`ene Bataillon\\ 34095 Montpellier CEDEX 5\\ France}
\email{burgunder@math.univ-montp2.fr} \urladdr{www.math.univ-montp2.fr/{$\sim$}burgunder/}

\title[A symmetric graph complex and Leibniz homology]{A symmetric version of Kontsevich graph complex and Leibniz homology} 
\subjclass[2000]{} 
\keywords{ Kontsevich graph complex, Leibniz homology, graph homology, Zinbiel-associative bialgebras, co-invariant theory}

\date{\today}
\begin{abstract} 
Kontsevich has proven that the Lie homology of the Lie algebra of symplectic vector fields can be computed in terms of  the homology of a graph complex. We prove that the Leibniz homology of this Lie algebra can be computed in terms of the homology  of a variant of the graph complex endowed with an action of the  symmetric groups. The resulting isomorphism is shown to be a Zinbiel-associative bialgebra isomorphism.  
\end{abstract}

\maketitle
\let\languagename\relax
In his papers \cite{K1} and \cite{K2}, Kontsevich proved that the homology of the Lie algebra $\Sp(\Com)$ of symplectic vector fields on a formal manifold
can be computed through graph homology: 
there exists a canonical co-commutative commutative bialgebra isomorphism 
\begin{equation*}\label{thm:K}
H_*(\Sp(\Com))\cong
\Lambda (H_*(\textrm{connected graph complex}))\ .  
\end{equation*}

In the literature, we can find another theorem of similar nature due to Loday, Quillen and  Tsygan cf.~\cite{LQ, T}.  It states that
the homology of the Lie algebra of matrices on an associative algebra can be computed as the exterior power of cyclic homology: there exists a canonical co-commutative commutative bialgebra isomorphism 
\begin{equation*}\label{thm:LQT}
H_*(\mathrm{gl}(A))\cong \Lambda (HC_{*-1}(A)) \ .
\end{equation*} 
 These two theorems are closely linked since the cyclic homology of an algebra can be seen as the graph homology of polygons labelled by elements of
the algebra.

Another similar theorem involving Leibniz homology $HL$ and Hochschild homology $HH$ has been proven by Cuvier and Loday (cf.~\cite{C, L}): there exists a vector space isomorphism 
$$ HL_*(\mathrm{gl}(A))\cong T (HH_{*-1}(A)) \ .$$
The Leibniz homology is the homology of the chain complex built over the tensor power $T(\gg)$, whereas the chain complex considered for the Lie homology is the exterior power $\Lambda(\gg)$ (quotient of $T(\gg)$ by the symmetric group action).  

The aim of this paper is to compute the Leibniz homology of the Lie algebra $\Sp(\Com)$. We construct a variant of Kontsevich graph complex, called the \emph{symmetric graph complex}. In dimension $n$ it is equipped with an action of the symmetric group $\Sigma_{n}$. Its quotient by this action gives Kontsevich graph complex. We show that there exists an  isomorphism:
\begin{equation}\label{eq:isomorphism}
   HL_*(\Sp(\Com))\cong T(H_*(\textrm{connected symmetric graph complex}))\ .
\end{equation}

On the left-hand side the direct sum of matrices induces an associative algebra structure, and the diagonal induces a Zinbiel coalgebra structure. They make $ HL_*(\Sp(\Com))$ into a Zinbiel-associative bialgebra. On the right-hand side there is an obvious free-cofree Zinbiel-associative bialgebra structure. Our isomorphism is shown to be an isomorphism of Zinbiel-associative bialgebras.

Under quotienting by the action of the symmetric groups our proof gives, as an immediate Corollary, a proof of Kontsevich theorem. Moreover, the two isomorphisms are related by a commutative diagram:

\begin{displaymath}
\xymatrix{
HL_*(\Sp(\Com))\ar[d]_{(\ )_{\Sigma_n}}&\cong &T(H_*(\textrm{connected symmetric graph complex}))\ar[d]^{(\ )_{\Sigma_n}}\\
H_*(\Sp(\Com))& \cong& \Lambda (H_*(\textrm{connected graph complex})) \ .}
\end{displaymath}

Kontsevich theorem leads to many types of generalisations. Hamilton and Lazarev proved it in the orthosymplectic context, cf.\ \cite{H},\cite{HL}. Mahajan extended it for reversible operads in \cite{M} and
 Conant and Vogtmann extended Konstevich proof to any cyclic operad, cf.\ \cite{CV}.  
We intend to show in a sequel to this paper, that such  generalisations are possible in the Leibniz homology context too. 

The paper is constructed as follows:  the first section sets the notations for Lie and Leibniz homology.  The second section introduces
the notion of symmetric graphs used in section three to state the main theorem.  The next sections are devoted to the proof of this
theorem.  Section five is the first step of the proof known as the Koszul trick, section six is devoted to recalls on co-invariant theory
for the symplectic algebra, that are then applied to the Leibniz complex of $\Sp(\Com)$.  The next section introduces chord diagrams with
a little digression to describe the chain complex of chord diagrams.  Section eight mimicks Kontsevich's idea to produce an isomorphism
between the classes of chord diagrams under a symmetric action and graphs.  But here the graphs that arise are labelled.  In section ten, we reduce the computation of the homology of
connected graphs to the homology of connected graphs with no bivalent vertices.
Then, in the last section we show that the isomorphism of the main theorem is an isomorphism of Zinbiel-associative bialgebras.
The appendix of this paper is devoted to a rigidity theorem, analogous to Hopf-Borel, for Zinbiel-associative bialgebras and its dual version.

In the sequel $\KK$ denotes a field of characteristic zero.

\section{Leibniz and Lie homology}\label{sectionLieLeibniz} In order to set the notations, we recall the basic notions of Leibniz and Lie
chain complexes, adjoint representation and co-invariants.

\begin{definition} A \emph{Leibniz algebra} $L$ is a vector space over a field $\KK$ endowed with a bilinear map $L\times L\rightarrow
L$, denoted $(x,y)\mapsto[x,y]$ and called the \emph{bracket} of $x$ and $y$, verifying the \emph{Leibniz identity}:  $$
[[x,y],z]=[[x,z],y]+[x,[y,z]] \textrm{ for all } x,y,z\in L\ .  $$ \end{definition}

This identity means that the operation $[-,z]$ is a derivation with respect to the bracket.

\begin{definition} A vector space $\gg$ over a field $\KK$, endowed with an operation $\gg\times\gg\rightarrow\gg$, denoted
$(x,y)\mapsto[x,y]$, is called a \emph{Lie algebra} over $\KK$ if $\gg$ is a Leibniz algebra and if the bracket verifies moreover that $$
[x,x]=0 \textrm{ for all } x\in\gg\ .  $$ \end{definition}

It is clear that the axiom $[x,x]=0$ for all $x\in\gg$ implies the anti-commutativity axiom, i.e.  $[x,y]=-[y,x]$ for all $x,y\in \gg$,
by applying the bilinearity hypothesis to the element $[x+y,x+y]$.  Moreover the identity satisfied by the Leibniz algebra implies the
\emph{Jacobi identity} :  \begin{equation*} [[x,y],z]+[[z,x],y]+[[y,z],x]=0 \textrm{ for all } x,y,z\in L\ , \end{equation*} under the
assumption of anti-commutativity.  \medskip

Let $\gg$ be a Lie algebra over $\KK$.  The Chevalley-Eilenberg chain complex is defined as follows:
$$\cdots\stackrel{d}{\longrightarrow} \Lambda^n \gg \stackrel{d}{\longrightarrow} \Lambda^{n-1} \gg \stackrel{d}{\longrightarrow}\cdots
\stackrel{d}{\longrightarrow} \Lambda^1 \gg\stackrel{d}{\longrightarrow} \KK\ ,$$ where $\Lambda^n \gg$ is the $n$th exterior power of
$\gg$ over $\KK$, and the map $d$ is given by the classical formula :  \begin{eqnarray*}d(g_1\wedge\cdots\wedge g_n):= \sum_{1\leq
i<j\leq n} (-1)^{i+j+1} [g_i,g_j]\wedge g_1\wedge\cdots\wedge \widehat{g_i}\wedge\cdots\wedge\widehat{g_j}\wedge\cdots\wedge g_n \ ,
\end{eqnarray*} where $\widehat{g_i}$ means that $g_i$ has been deleted.  The homology of the Chevalley-Eilenberg chain complex is
denoted $H_*(\gg,\KK)$ or $H_*(\gg)$ if no confusion occurs with this notation.

Let $\gg$ be a Leibniz algebra over $\KK$.  The Leibniz chain complex $\mathrm{CL}_*(\gg)$ has been defined by J.-L.  Loday in
\cite{Lleib} as a lifting of the Chevalley-Eilenberg complex by:  $$\cdots\stackrel{d}{\longrightarrow} \gg^{\t n}
\stackrel{d}{\longrightarrow} \gg^{\t n-1} \stackrel{d}{\longrightarrow}\cdots \stackrel{d}{\longrightarrow}\gg
\stackrel{d}{\longrightarrow}\KK$$ where $\gg^{\t n}$ is the $n$th tensor power of $\gg$ over $\KK$, and where the map $d$ is given by
the following formula :  \begin{eqnarray*} d(g_1\t\cdots\t g_n):=\sum_{1\leq i<j\leq n} (-1)^{j} g_1\t\cdots\t g_{i-1}\t[g_i,g_j]\t
g_{i+1}\t\cdots\t\widehat{g_j}\t\cdots\t g_n \ .  \end{eqnarray*} In this paper we denote the homology of the Leibniz chain complex by
$HL_*(\gg,\KK)$ or $HL_*(\gg)$ if it doesn't lead to any confusion.

\section{Symmetric graph homology}

We introduce the notion of symmetric graph complex which gives rise to  symmetric graph
homology. The graphs that we consider are labelled by integers and endowed with an action of the symmetric groups on the labellings. We show moreover that the chain complex and the homology of graphs admit a structure of generalised bialgebra, namely a Zinbiel-associative bialgebra structure.

\begin{definition}\label{def:graph}
 Let $m$ be an integer. An \emph{oriented symmetric graph} $G$ is a triplet made of the ordered set $V(G):=\{1,\ldots,m\}$ of vertices, a
set of edges $E(G)$ and a map $\alpha_G:E(G)\to V(G)\times V(G)$.  The symmetric graph $G$ is said to be \emph{non-oriented} if
$\alpha_G:E(G)\to S^2(V(G))$.

To ease the writing the set $V(G)$ is referred to the set of \emph{vertices} of $G$, and the set $E(G)$ is referred to the set of
\emph{edges} of $G$.  A vertex $v \in V(G)$ has \emph{valency} $n$ if the cardinality of the set $\{e\in E(G)|\exists a\in V(G) :
\alpha_G(e)=(a,v) \textrm{ or } (v,a)\}$ is $n$.  The elements of $\alpha_G(e)$ are called the \emph{incident half-edges} of $e$.  An
edge whose composite half-edges are incident to the same vertex is called a loop.

A \emph{labelled graph} is a graph $G$ together with a map from the set of vertices $V(G)$ to the set of labellings $\{p_1,q_1,\ldots\}$.

The symmetric group $\Sigma_m$ acts on the graph $G=(\{1,\ldots,m\},\{e_i\},\alpha_G)$ as follows :
let $\alpha_G(e_i)=(i_1,i_2)$
$$\sigma\cdot(V(G),E(G),\alpha_G):=\sgn(\sigma)( V(G),E(G),\sigma\cdot \alpha_G)\ ,$$
where $
(\sigma\cdot \alpha_G)(e_i)=(\sigma(i_1),\sigma(i_2))\ .
$
\end{definition}

\begin{definition} Two graphs $G_1$ and $G_2$ are said to be isomorphic if :  
\begin{enumerate} 
\item $V(G_1)=V(G_2)$, 
\item $|E(G_1)|=|E(G_2)|$, 
\item $\Im \alpha_{G_1}=\Im \alpha_{G_2}$.  
\end{enumerate} 
\end{definition}

In the sequel, we will consider isomorphic classes of non-oriented symmetric graphs without loops, and isomorphic classes of oriented
symmetric graphs without loops (unless otherwise stated).  We will refer to them as graphs and oriented graphs respectively.

\begin{notation} The set of all finite graphs will be denoted $\GG$, and the set of finite oriented graphs is denoted $\stackrel{\to}{\GG}$.  The set of finite graphs with $m$ vertices (i.e.  every graph $G$ such that $|V(G)|=m$) is denoted $\GG_m$.  The set of graphs such that the valence of their vertices is respecctively  $k_1,\ldots,k_m$ is denoted
$\GG_{k_1,\ldots,k_m}$.  
\end{notation}

\begin{example}\label{ex:graphexample} The connected graph $G=(V(G),E(G),\alpha_G)$ with $V(G)=\{1,2,3\}$, $E(G)=\{e_1,\ldots, e_4\}$ and
$\alpha_G(e_1)=\{1,3\},\alpha_G(e_2)=\{1,2\},\alpha_G(e_3)=\{1,2\},\alpha_G(e_4)=\{2,3\}$ can be represented geometrically as in figure
\ref{figure:graphexample}.  \begin{figure}[htbp] \begin{center} \input{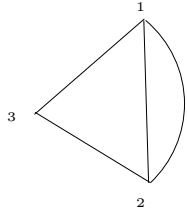} \caption{Geometric interpretation of the
graph $G$} \label{figure:graphexample} \end{center} \end{figure} \end{example}

\begin{definition}\label{def:operationgraph} Let $G_1$ and $G_2$ be two graphs with $V(G_i)=\{1,\ldots,n_i\}$, $E(G_i)=\{e_j\}_{1\leq j\leq r_i}$.  The ordered
disjoint union of two graphs $G_1\cdot G_2$ is defined as follows: 
\begin{eqnarray*} 
V(G_1\cdot G_2)&:=&\{1,\ldots,n_1+n_2\}\ ,\\
E(G_1\cdot G_2)&:=&\{e'_j\}_{1\leq j\leq r_1+r_2}\ , \\ 
\alpha_{G_1\cdot G_2}(e'_j)&:=&
\left\{\begin{array}{ccc}
\alpha_{G_1}(e_j)&\textrm{if } j\leq r_1\ ,\\
\alpha_{G_2}(e_{j-r_1})&\textrm{if }j\geq r_1+1 \ .  
\end{array}\right . 
\end{eqnarray*}
This operation endows the vector space $\KK[\GG]$ with a structure of associative (non-commutative) algebra.

A graph $G$ is \emph{connected} if for any graph $U,V\in\GG$, $U\cdot V\neq G$.  We denote $\GG_c$ the set of connected graphs.  We will
denote $\GG_c^3$ the set of connected graphs such that every vertex is of valence at least $3$.

It is to be noted that in our definition there is a notion of order on the connected components of a graph.  Indeed 
let $G_1$ and $G_2$ be two components of a graph $G$, then $G_1$ is greater than $G_2$ if the minimum of the labels of the  vertices in $G_1$ is greater than the minimum of the labels of the vertices
in $G_2$.  This gives rise to a vector space isomorphism between the tensor module $T(\KK[\GG_{c}])$ and $\KK[\GG]$.
\end{definition}

The definition of graphs of definition \ref{def:graph} encodes non-necessarily connected graphs, as illustrated in the following example.
  \begin{example}\label{ex:nonconnexample} Let
$H=(V(H),E(H),\alpha_H)$ be the graph defined as $V(H)=\{1,2,3,4\}$, $E(H)=\{e_1,\ldots,e_4\}$ and
$\alpha_H(e_1)=\{1,2\},\alpha_H(e_2)=\{1,2\},\alpha_H(e_3)=\{3,4\},\alpha_H(e_4)=\{4,3\}$.  This graph can be represented geometrically
as in figure \ref{figure:nonconectgraph}, taking into account the order of the components.  \begin{figure}[htbp] \begin{center}
\input{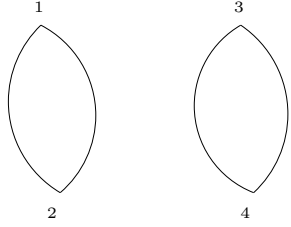} \caption{Geometric interpretation of the graph $H$} \label{figure:nonconectgraph} \end{center}
\end{figure}

This graph is the disjoint union of $H_1\cdot H_1$, where $H_1$ is the graph defined by $V(H_1):=\{1,2\}$, $E(H_1)=\{e_1,e_2\}$ and
$\alpha_{H_1}$ satisfies :  \begin{eqnarray*} \alpha_{H_1}(e_1)=\{1,2\}\quad \alpha_{H_1}(e_2)=\{1,2\}\ .  \end{eqnarray*}

Since $H$ is a class of isomorphic symmetric graphs, we can write $H$ as follows $V(H)=\{1,2,3,4\}$ and
$\Im\alpha_H=\{\{1,2\},\{1,2\},\{3,4\},\{4,3\}\}$ for sake of simplicity.  Therefore with this notation $H_1$ is defined as
$(\{1,2\},\{\{1,2\},\{1,2\}\})$.  \end{example}

\begin{proposition} The associative algebra $(\KK[\GG],\cdot)$ is isomorphic to the free associative algebra $T(\KK[\GG_c])$.
\end{proposition}

\begin{proof} Any non-connected graph can be seen as the ordered disjoint union of some uniquely determined non-empty connected graphs.
The ordered disjoint union is an associative non-commutative operation on graphs.  This gives rise to an isomorphism between the tensor
algebra over the connected graphs and the algebra of graphs.  \end{proof}

In order to define the differential in the graph complex, we first describe how to contract an edge in a given graph.  \begin{definition}
Let $G=(\{1,\ldots,m\},\{e_j\}_{1\leq j\leq n},\alpha_G)$ be a graph and let $e=(i,j)$ be one of its edges, which we assume not being a
loop.  We define a new graph $G/e$ as follows:  \begin{enumerate} \item the set of edges of $G/e$ are the edges of $G$ minus $e$, \item
the set of vertices of $G/e$ is $\{1,\ldots, m-1\}$, \item the map $\alpha_{G/e}:E(G/e)\to V(G/e)\times V(G/e)$ is defined as follows :
$$\alpha_{G/e}(e')=(\std\t\std)\circ\alpha_{G}(e')$$ where, \begin{equation*} \std(k)=\left\{\begin{array}{ccc} k&\textrm{if }
k<\min(i,j)\ ,\\ k-1&\textrm{if }i<k \textrm{ and }k \neq \max(i,j)\\ \min(i,j)&\textrm{if }k=\max(i,j)\ .  \end{array}\right .
\end{equation*} 
\end{enumerate} The map $\std$ is known as the standardisation map.  \end{definition} This definition do not depend on the
representative $G$ in the class of isomorphic graphs.

\begin{example} We consider the graph $G$, defined by $V(G)=\{1,\ldots,4\}$ and $\Im\alpha_{G}=\{\{1,4\},\{1,3\},\{2,4\},\{2,3\},\{3,4\}\}$, and
the edge $e=\{1,4\}$.  The resulting graph $G/e$ is defined as $V(G)=\{1,2,3\}$ with $\Im\alpha_{G}=\{\{1,3\},\{2,1\},\{2,3\},\{3,1\}\}$.

\begin{figure}[htbp]
 \begin{center}
  \input{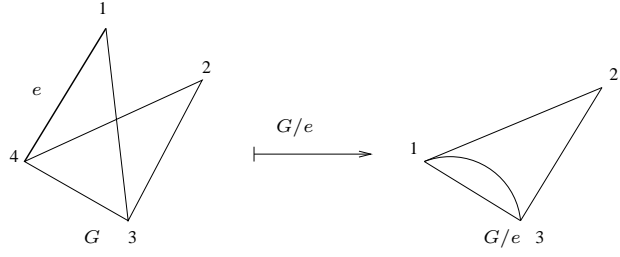}
   \caption{Contracting vertex $e$ in graph $G$.}
     \label{figure:smashgraph}
  \end{center}
 \end{figure}
\end{example}

\begin{proposition}\label{prop:equivalence} Let $G\in\GGto$ be a graph obtained by $G'\in\GGto$ by a change of orientation on a edge $e=[i,j]$ into $[j,i]$. We set $G\sim -G'$.  

The algebra spanned by the oriented graphs quotiented by the above equivalence relation 
is exactly the algebra spanned by the graphs without loops.  \end{proposition} 
\begin{proof} It is clear that
the loops are null.  Moreover the orientation is mod out by the sign relation on orientation.  \end{proof}

To ease the notation we will define a sign $\epsilon(e)$ for any oriented edge $e=[i,j]$ which will depend on the orientation of the edge and the
number associated to the vertex of each incident half-edge.  
\begin{definition} 
Let $G\in\GGto$ be an oriented graph, and let $e=[i,j]$ be an
oriented edge.  We define the sign $\epsilon(e)\in\{-1,1\}$ as follows:  
\begin{equation*} 
\epsilon(e)=
\left\{\begin{array}{ccc}
-1&\textrm{if } i>j\ ,\\
 1&\textrm{if }i<j\ .
\end{array}
\right .
\end{equation*}
\end{definition} 

We are now able to describe the
complex of oriented symmetric graphs $C(\GG,\delta)$.

\begin{definition}\label{def:cooperationgraph}
The $n$th chain module $C_n(\GGto)$ is $C_n(\GGto):=(\KK[\GGto])^{\t n}$.  The differential on the chain module is
defined as follows:   let $\stackrel{\to}{G}$ be an oriented graph,
$$\delta( \stackrel{\to} {G}):=\sum_{e}(-1)^{\textrm{max}(i,j)} \epsilon(e)  \stackrel{\to}{ G/e}\ ,$$ where the sum runs
over all edges $e=[i,j]$.
\end{definition}

\begin{proposition}
The chain complex $C_*(\stackrel{\to}{\GG})$ passes through the quotient by the equivalence relation defined in proposition \ref{prop:equivalence}, defining the chain complex associated to non-oriented graphs . 
\end{proposition}

\begin{proof} 
The differential does not depend on the oriented representative of the graph $G$.  Indeed, let $ \stackrel{\to} {G}$ be an oriented representative of the graph $G$.  Let $v$ be a vertex in $ \stackrel{\to} {G}$ with orientation $[i,j]$.  Consider the graph $\stackrel{\to}{ G'}$ which is the same graph as $ \stackrel{\to} {G}$ with the orientation of $v$ being $[j,i]$.  Therefore, direct
computation gives :  
\begin{eqnarray*} 
\delta( \stackrel{\to} {G})-\delta( \stackrel{\to} {G'})&=&
(-1)^{\textrm{max}(i,j)} (\epsilon(v) [\stackrel{\to} {G/v}]
-\epsilon(v') [ \stackrel{\to} {G'/v'}])\\
 &=& (-1)^{\textrm{max}(i,j)}\epsilon(v)(1-1)[ \stackrel{\to} {G/v}]\\
&=&0\ . 
 \end{eqnarray*} 
This end the proof as for more than one change of orientation, we will have just a sum of the above equality.
\end{proof}

\begin{example}\label{ex:graphdiff} 
We illustrate the differential of a graph by performing the calculation on the connected graph
defined in \ref{ex:graphexample}.  It gives the following sum $-2(\{1,2\},\{\{1,2\},\{1,2\},\{1,2\}\})$.  The sum can be reinterpreted
geometrically as in figure \ref{figure:graphdiff} where we have just contracted each edge \begin{figure}[htbp] \begin{center}
\input{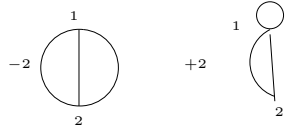} \caption{Sum of contracted edges of $G$} \label{figure:graphdiff}\end{center} \end{figure} but we have also to
take into account that the loops are null.  Figure \ref{figure:graphdiffresult} gives the geometrical result.  \begin{figure}[htbp]
\begin{center} \input{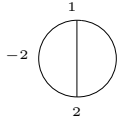} \caption{The differential of the graph $G$} \label{figure:graphdiffresult}\end{center}
\end{figure}

\end{example}

Our aim is to prove that the graph complex admits a structure of Zinbiel-associative bialgebra.  The definition and the main properties
of Zinbiel-associative bialgebras can be find in the appendix.  It includes a structure theorem stating that a connected
Zinbiel-associative bialgebra can be reconstructed from its primitives.  This theorem is analogous to the Milnor-Moore theorem for
co-commutative commutative bialgebras.

First, we define the Zinbiel coalgebra structure on the graph complex as the canonical Zinbiel coproduct on the tensor module $T(\KK[\GG_{c}])$.  
\begin{definition}\label{def:coop_graph} 
Let $\DD:C(\GG,\delta)\rightarrow C(\GG,\delta)\t C(\GG,\delta)$ be defined as the co-half shuffle on the ordered disjoint union of graphs.  That is to say, let
$G=G_1\ldots G_n$ be the ordered disjoint union of the connected graphs $G_i$, for $1\leq i \leq n$, the co-half shuffle is defined as:
$$\DD(G):=G_1\SS(G_2\ldots G_n)=G_1\sum_{p+q=n}\sum_{\underline i \in Sh_{p,q}}G_{i_1}\ldots G_{i_p}\t G_{i_{p+1}}\ldots G_{i_n}$$
where the sum is extended over all $(p,q)$-shuffles $\underline i$ (i.e.  in the multi-indices $\underline i=(i_1,\ldots,i_n)$ the integers $2,
\ldots, p$ are ordered and so are $p+1, \ldots, n$).  
\end{definition}

\begin{example} 
Let $G$ be the graph defined by
$V(G)=\{1,\ldots,8\}$ and $\Im\alpha_{G}=\{\{7,8\}\{1,3\},\{4,5\},\{1,3\},\{1,2\},\{4,5\}\{2,3\},\{6,8\},\{7,8\},\{6,7\},\{6,8\}\}$.  It is clear
that $G=G_1\cdot G_2\cdot G_3$.  Indeed, $G_1=(\{1,2,3\},\{\{1,3\},\{1,3\},\{1,2\},\{2,3\}\})$, $G_2=(\{1,2\},\{\{1,2\},\{1,2\}\})$ and
$G_3=(\{1,2,3\},\{\{1,3\},\{1,3\},\{1,2\},\{2,3\},\{1,3\}\})$.
Applying the coproduct to this graph gives the following:  
$$G_1\cdot G_2\cdot G_3\t 1+G_1\t G_2\cdot G_3+G_1\cdot G_2\t G_3
+G_1\cdot G_3\t G_2\ .  $$
This can be geometrically interpreted as in figure \ref{figure:coproduitgrap}.

\begin{figure}[htbp] \begin{center} \input{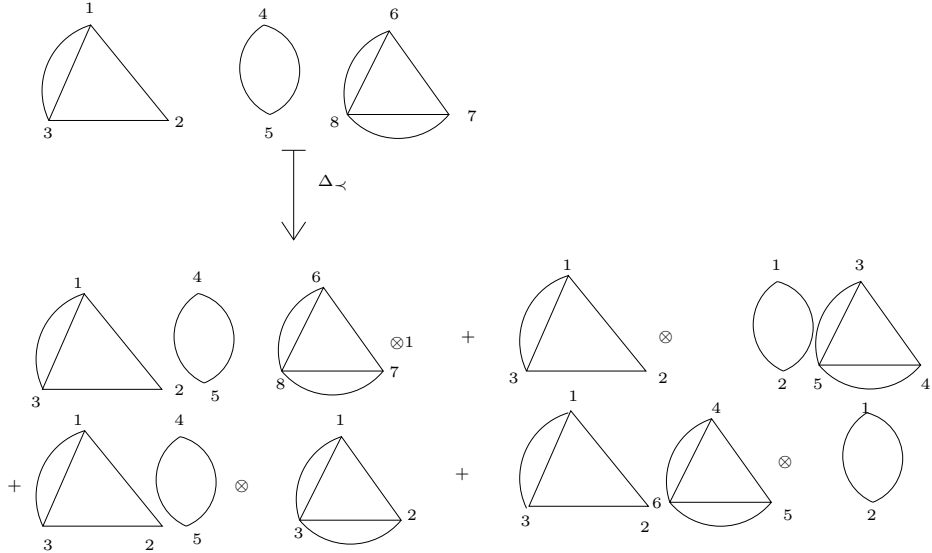} \caption{Coproduct of the graph $G$}
\label{figure:coproduitgrap}\end{center} \end{figure}

\end{example}

The Zinbiel coalgebra can also be endowed with an associative algebra structure.  This endows the chain complex with a
Zinbiel-associative bialgebra, as shown in the following proposition.  

\begin{proposition} 
Let $\cdot:C_p(\GG)\t C_q(\GG)\to
C_{p+q}(\GG)$ be the disjoint ordered union of graphs.  This operation endows the chain complex with an associative algebra structure.
Together with the coalgebra structure, defined in definition \ref{def:coop_graph}, it endows the chain complex of symmetric graphs with a structure of
Zinbiel-associative bialgebra.  
\end{proposition}
 
\begin{proof}
 The operation $\cdot $ is clearly associative (and not commutative).  The
co-half shuffle is a Zinbiel coproduct cf. appendix \ref{App:As-Zinb}.  It suffices to verify the compatibility relation.  This can be
done by dualising the arguments of \cite {B}, see appendix \ref{App:Zinb-As}.  
\end{proof}

\begin{proposition}\label{prop:primgraph}
 The primitive part of $C(\GG,\delta,\DD)$ is the subvector space spanned by connected graphs.
\end{proposition}

\begin{proof} It's known that the primitive part of $T(V)$ for the co-half shuffle is $V$ (cf.\ appendix \ref{App:Zinb-As} for a proof).  This remark
ends the proof.  
\end{proof}

\begin{corollary} 
The homology $H(\GG,\cdot,\delta,\DD)$ of the graph complex is isomorphic to the graded Zinb-As bialgebra
$T(H(\GG_{\rm{c}}))$.  
\end{corollary} 

\begin{proof} 
It suffices to apply corollary \ref{cor:homandprim} of appendix \ref{App:Leib-Zinb} and proposition
\ref{prop:primgraph}.  
\end{proof}

\section{A Kontsevich analogue for the Leibniz homology} 

In this section, we introduce the Lie algebra $\Sp(\Com)$ and state the main theorem:  The Leibniz homology of the Lie algebra
$\Sp(\Com)$ can be computed thanks to the homology of the symmetric graph complex defined in the previous section.

\begin{definition} 
Let $m$ be a positive integer.  
Let $V_m$ be the $\KK$-vector space generated by $2m$ indeterminates
$\{p_1,\ldots,p_m,q_1,\ldots,q_m\}$ endowed with the standard symplectic form $\omega=\sum_{i=1}^{m}\mathrm{d}p_i\wedge \mathrm{d}q_i$.
We define the Lie algebra $\Sp_m(\Com)$ as follows:  the underlying vector space is $S^{\geq 2}(V_m):=\oplus_{k\geq 2}(\Vt k_m)_{S_k}$,
the module of commutative polynomials in indeterminates $p_1,\ldots,p_m,q_1,\ldots,q_m$.  It is endowed with the canonical Poisson bracket:
\begin{equation*} \{F,G\} = \sum_{i=1}^{n} \frac{\partial F}{\partial p_i}\frac{\partial G}{\partial q_i}-\frac{\partial G}{\partial
p_i}\frac{\partial F}{\partial q_i} \ , \end{equation*} where $F$ and $G$ are polynomials in indeterminates
$p_1,\ldots,p_m,q_1,\ldots,q_m$.

The inductive limit of these Lie algebras is denoted $\Sp(\Com):=\cup_{m\geq 1}\Sp_m(\Com)$.  \end{definition} Conant and Vogtman defined
in \cite{CV} a functor from cyclic operads to symplectic Lie algebras.  For $\PP=\Com$ it turns out to be the Lie algebra described
above.  In a future article, we will explain the choice of the notation $\Sp(\Com)$.

The algebra $\Sp_m(\Com)$ is a Lie algebra.  It is in particular a Leibniz algebra, and therefore we can consider its Leibniz homology.

Our aim is to prove the following theorem:

\begin{theorem}\label{thm:LeibnizKont}Let $\KK$ be a characteristic zero field.  There exists a canonical Zinbiel-associative bialgebra
isomorphism :  $$HL_*(\Sp(\Com))\cong T(H_*( \GG^3_c))$$ \end{theorem}

\begin{proof} The proof will be decomposed into four steps as follows :

First step:  quotient the chain complex $\mathrm{CL}_*(\Sp(\Com))$ by the action of the Lie reductive algebra $\Sp_{2m}(\KK)$.  This step is
known as the Koszul trick.

Second step:  apply the co-invariant theory and then mimic Kontsevich's idea to reduce the chain complex to a complex of graphs.

Third step:  reduce the chain complex of graphs to a smaller one by constructing explicit homotopies.

Fourth step: show that the isomorphism obtained is a Zinbiel-associative bialgebra isomorphism.  \end{proof}

\section{First step:  Koszul trick} The Leibniz chain complex of the Lie algebra $\Sp(\Com)$ can be reduced by the Koszul trick.  It is quasi-isomorphic to the chain complex spanned by the symplectic co-invariants.

\subsection{Adjoint representation and Leibniz homology} In this section we translate some well known homological properties of  Lie
algebras in the Leibniz context.  Most of this subsection can be found in \cite{L}.

\begin{proposition}\label{trivial} Let $\gg$ be a Leibniz algebra.  The adjoint action of $\gg$ on itself given by \begin{equation*}
[g_1\t\cdots\t g_n, g ]:=\sum_{i=1}^{n} g_1 \t\cdots\t [ g_i,g ]\t\cdots\t g_n \ , \end{equation*} is compatible with $d$.  The induced
action on $HL_*(\gg,\KK)$ is trivial.  \end{proposition} \proof The compatibility of the adjoint action with the differential is proved
in \cite{L} lemma 10.6.3 (10.6.3.0) by induction :  $$ d[\alpha,g]=[d\alpha,g] \textrm{ for all } \alpha\in\gg^{\t n} \textrm{ and }
g\in\gg \ .  $$ In order to prove the second assertion, we construct a homotopy $\sigma$ from the adjoint action to zero.

For any $y\in\gg$ let $\sigma(y):\gg^{\t n}\longrightarrow \gg^{\t n+1}$ be the map of degree one given by:  \begin{equation*}
\sigma(y)(\alpha):=(-1)^{n}\alpha\t y, \ \ \alpha\in \gg^{\t n} \ .  \end{equation*} Then, the following equality holds :
\begin{eqnarray*} d\sigma(y)(\alpha)+\sigma(y) d(\alpha)&=& \sum_{i,j=1}^{n+1} (-1)^{j+n+1} \alpha_1\t\cdots \t[
\alpha_i,\alpha_j]\t\cdots\t\widehat{\alpha_j}\t\cdots \alpha_n \t y\\ &&+\sum_{i,j=1}^{n} (-1)^{j+n} \alpha_1\t\cdots
\t[\alpha_i,\alpha_j]\t\cdots\t\widehat{\alpha_j}\t\cdots \alpha_n \t y\\ &=& [\alpha,y] \ .  \end{eqnarray*} This proves that
$\sigma(y)$ is a homotopy from $[-,y]$ to $0$, whence the assertion.  $\Box$

\begin{definition} Let $G$ be a group and $V$ be a left $\KK[G]$-module.  By definition the vector space of \emph{invariants} is
$$V^G:=\{v\in V|g\cdot v=v \ \textrm{for all } g\in G\}\ ,$$ and the space of co-invariants is :  $$V_G:=V/_{\{g\cdot v-v \}} \ .$$
\end{definition} When $G$ is finite and $|G|$ is invertible in $\KK$, there is a canonical isomorphism $V_G\cong V^G$ given by the
averaging map :  $$[v]\mapsto \frac{1}{|G|}\sum_{g\in G}g\cdot v \ .$$

\begin{proposition}\label{prop:quasiiso} Let $\gg$ be a Leibniz algebra.  Let $\hh$ be a reductive sub-Lie algebra of $\gg$.  Then the
surjective map $\gg^{\t n} \rightarrow (\gg^{\t n})_\hh$ induces an isomorphism on homology, \begin{equation*} HL_*(\gg,\KK)\cong
HL_*((\gg^{\t n})_\hh,d) \ .  \end{equation*} \end{proposition} \proof Since $\gg$ is a completely reducible $\hh$-module, the module
$\gg^{\t n}$ splits, as a representation of $\hh$, into a direct sum of isotypic components.  The component corresponding to the trivial
representation is $ (\gg^{\t n})_\hh$ (co-invariant module).  Let us denote $L_n$ the sum of all the other components.  Since $\hh$ is a
sub-Lie algebra of $\gg$ it is in particular a Leibniz algebra and then proposition \ref{trivial} implies that $d$ is compatible with the
action of $\hh$.  So there is a direct sum decomposition of complexes:  \begin{equation*} T(\gg)\cong ( T(\gg))_\hh\oplus L_\star.
\end{equation*} To finish the proof it suffices to prove that $L_\star$ is an acyclic complex.

Since $L_\star$ is made of simple modules which are not trivial $\hh$-modules, the components of $HL_*(L_*)$ are not trivial either.  But
by proposition \ref{trivial}, $HL_*(\gg,\KK)$ is a trivial $\gg$-module and so a trivial $\hh$-module.  Therefore $HL_*(L_\star)$ has to
be zero.$\Box$

The same properties are true for Lie algebras with the Chevalley-Eilenberg chain complex, see \cite{L} section 10.1.8.  \subsection{First
step of the proof} We apply the above theory to the symplectic Lie algebra $\Sp(\Com)$ with the action of the reductive symplectic Lie
algebra $\Sp(\KK)$.  First we verify that the symplectic Lie algebra $\Sp(\KK)$ is included in the Lie algebra $\Sp(\Com)$, cf \cite{CV}.

Recall that $\Sp(2n)$ is the set of $2n \times 2n$-matrices $a$ satisfying $aj+j{}^ta=0$ where $j= \left( \begin{array}{cc}
0 & \Id \\ -\Id & 0 \end{array}\right)$.

\begin{proposition}\label{prop:sp2}\cite{CV} The symplectic Lie algebra $S^2(V_m)$ is isomorphic to the Lie algebra $\Sp_{2n}(\KK)$.
\end{proposition} 

\commentc{\begin{proof}This proof can be found in \cite{CV}.  Note that $S^2(V)$ is a symplectic Lie algebra as it is
stable under the bracket.  Indeed, the bracketing of two elements of $S^2(V)$ is a sum of at most two elements in $S^2(V)$.

Recall that $\Sp(2n)$ is the set of $2n \times 2n$-matrices $a$ satisfying $aj+j{}^ta=0$ where $j={\mathcal A} = \left( \begin{array}{cc}
0 & \Id \\ -\Id & 0 \end{array}\right)$.  Define $E_{ij}$ to be the $n\times n$-matrice with a one in the $(i,j)$-entry and zero
elsewhere.  The following map $S^2(V)\to \rm{End} (V):vw\mapsto \oo(v,-)w+\oo(w,-)v$ is a Lie algebra isomorphism.  Indeed its sends the
basis of $S^2(V)$ on the basis of $\Sp_{2n}(\KK)$:  \begin{eqnarray*}&&p_iq_j\mapsto \left( \begin{array}{cc} -E_{ij} &0 \\ 0 & E_{ji}
\end{array}\right)\ ,\\ &&p_ip_j\mapsto \left( \begin{array}{cc} 0 &0 \\ E_{ij}+E_{ji} &0 \end{array}\right)\ ,\\ &&q_iq_j\mapsto \left(
\begin{array}{cc} 0 &-E_{ij}-E_{ji} \\ 0 &0 \end{array}\right)\ .  \end{eqnarray*} Moreover it is stable under the bracketing.
\end{proof}}

Since $\Sp_{{2m}}(\KK)$ is contained into $\Sp_{2m}(\Com)$ it acts on $\Sp_{2m}(\Com)$ by the inner action.

 \begin{corollary} The Leibniz complex of the Lie algebra $\Sp(\Com)$ is quasi-isomorphic to the Leibniz complex
$\Sp(\Com)_{\Sp(\KK)}$.  \end{corollary} 

\begin{proof} Apply proposition \ref{prop:quasiiso} to the Leibniz algebra $\gg=\Sp_{2m}(\Com)$, and
the reducible algebra $\hh=\Sp_{2m}$ which is included in $\Sp_{2m}(\Com)$ by the above proposition.  Then, take the inductive limit.
\end{proof}

\section{The co-invariant theory for the symplectic group \cite{L,GW,P}}

 Co-invariant theory for the symplectic group is the main key to
reduce our complex to a complex of graphs.  The first part of this section is devoted to some recalls on the co-invariant theory for the
symplectic group, due to Procesi \cite{P}.  Then, we apply this theory to our Leibniz chain complex.  

\subsection{Recall on co-invariant theory for the symplectic group} 

\begin{definition} 
Let $V$ be a vector space over $\KK$ and let $\oo:V\times V\rightarrow \KK$ be a non-degenerate
skew-symmetric bilinear form.  We define the \emph{symplectic group} relative to $\oo$ as a subgroup of the general linear group on $V$,
denoted $GL(V,\KK)$, by:  $$ Sp(V,\oo):=\{g\in GL(V,\KK):\oo(gx,gy)=\oo(x,y), \textrm{ for all }x,y\in V\} \ , $$ 
\end{definition}

\begin{definition} 
Consider the polynomial algebra $\KK[y_{ij}]$, where $i\neq j$ ranges over $\{1,\ldots, 2n\}$, and where the relation
$y_{ij}=-y_{ji}$ holds.  Since we are in characteristic zero, it implies $y_{ii}=0$.  Define $A_n$ as the subspace spanned by the
monomials $y_{i_1 i_2}\ldots y_{i_{2n-1} i_{2n}}$ such that $\{i_1,\ldots,i_{2n}\}$ is a permutation of $\{1,\ldots,2n\}$.
\end{definition}

\begin{definition} 
The symmetric group $\Sigma_n$ acts on the left on $A_n$ by permuting the variables as follows:  $$ \sigma\cdot y_{i_1
i_2}\ldots y_{i_{2n-1} i_{2n}}= y_{\sigma^{-1}(i_1) \sigma^{-1}(i_2)}\ldots y_{\sigma^{-1}(i_{2n-1}) \sigma^{-1}(i_{2n})}\ , $$ for all
$y_{i_1 i_2}\ldots y_{i_{2n-1} i_{2n}}\in A_n$ and for all $\sigma\in\Sigma_n$.  
\end{definition} 

Dualising the assertions in \cite{L}
section 9.5, leads to the following formulations of the two fundamental theorems of co-invariant theory in the symplectic context.

\begin{theorem}[First Fundamental Theorem for the Symplectic group]\label{theo:FFT}
Let $V$ be a finite-dimensional vector space over
$\KK$ and $\oo:V\times V\rightarrow \KK$ be a non-degenerate skew-symmetric bilinear form.  The map $T^*:(\Vt {2r})_{\mathit{Sp}
(V)}\longrightarrow A_r$ induced by :  $$v_1\t\cdots\t v_{2r}\mapsto \sum_{y_{i_1 i_2}\cdots y_{i_{2r-1}i_{2r}}\in A_r}\oo^{\t
r}(v_{i_1}\t v_{i_2}\t\cdots\t v_{i_{2r-1}}\t v_{i_{2r}})\ y_{i_1 i_2}\cdots y_{i_{2r-1}i_{2r}}\ ,$$ where the sum is over all monomials
of $A_r$, is injective.  
\end{theorem}

\begin{theorem} [Second Fundamental Theorem for the Symplectic Group]\label{theo:SFT}
Let $V$ be a finite-dimensional vector space over
$\KK$ and $\oo:V\times V\rightarrow \KK$ be a non-degenerate skew-symmetric bilinear form.  The co-kernel of the map $T^*:(\Vt
{2r})_{\mathit{Sp} (V)}\longrightarrow A_r$ is the sum of all the isotypic components $((A_r)_\lambda)^*$,
$\lambda=\{\lambda_1,\cdots,\lambda_r\}$ such that $\lambda_1\geq \frac{1}{2} \dim V +2$.  In particular $T^*$ is an isomorphism as soon
as $\dim V\geq 2r$.  
\end{theorem}

These two theorems can be summarized in the following short exact sequence:  $$ 0 \longrightarrow (\Vtm
{2r})_{\mathit{Sp}_{2m}(\KK)}\longrightarrow A_r\longrightarrow \oplus_\lambda ((A_r)_\lambda )^*\longrightarrow 0 \ .$$

\begin{proposition}\label{odd}
Let $V$ be a finite-dimensional vector space over $\KK$ of characteristic zero and $\oo:V\times
 V\rightarrow \KK$ be a non-degenerate skew-symmetric bilinear form.  If n is odd, then we have $(\Vt n)_{\mathit{Sp}(V)}=0$.
 \end{proposition}

\begin{proof}This proof can be found in \cite{GW}. For $(V,\rho)$ a representation, then define $\rho_n$ as the representation of the tensor product $V^{\t n}$.  Clearly,
there are no invariants:  as $-I\in \Sp(V)$, then for all $x\in \Vt n$ we have $\rho_n(-I)(x)=(-1)^n x=-x$, since $n$ is odd.  As for
finite dimensional vector spaces over a characteristic zero field, invariants and co-invariants are isomorphic, the assertion is true.
\end{proof}

Let $M$ be $\mathrm{Sp}(\KK)$-bimodule, there is a canonical structure of $\Sp(\KK)$-module on $M$.  Moreover, we have the following
isomorphism of co-invariant spaces $M_{\mathit{Sp}(\KK)}=M_{\Sp(\KK)}$.

\subsection{Co-invariant theory applied to the Leibniz complex of $\Sp(\Com)$} We will now apply the above theory of co-invariants for
 the symplectic group to reduce the Leibniz chain complex of $(\Sp(\Com))_{\Sp(\KK)}$ to a more explicit chain complex.  As we focus on
 the inductive limit of $\Sp(\Com)$ the two fundamental co-invariant theorems give rise to an isomorphism.

If we focus on the non-stable case, that is to say $\Sp_m(\Com)$ where $m$ is fixed, then the co-invariant theorems would only stabilise
one part of each chain module.  It is very different from the Loday-Quillen case where they compute the homology of the general linear group of an associative algebra.  In their case the co-invariant theory gives rise to an isomorphism $H_{n}(\mathrm{gl}_{n}(A))\cong H_{n+1}(\mathrm{gl}_{n}(A))\cong\ldots$.  Therefore, they can compute the first obstruction to stability.

In order to apply the co-invariant theory, we have to verify first that the action of $\Sp(\KK)$ and of the symmetric group commute.

\begin{proposition}\label{actioncommute} Let V be a finite dimensional vector space.  Let $m=\sum_{i=1}^n k_i$.  The actions of $\Sp(V)$
on $\Vt m$ and of $\Sigma_{k_1}\times\cdots\times\Sigma_{k_n}$ on $\Vt m$ commute.  \end{proposition} \proof Let $A\in \Sp(V)$,
$x:=x_1\t\cdots\t x_n\in \Vt n$, $\sigma\in\Sigma_{k_1}\times\cdots\times\Sigma_{k_n}$.  We denote $\rho$ the action of the symplectic
Lie algebra.  We show by direct computation that the two actions commute.  \begin{eqnarray*} \rho(A)\circ\sigma x&=&\rho(A)\cdot
x_{\sigma^{-1}(1)}\t\cdots\t x_{\sigma^{-1}(n)}\\ &=& \sum_{j=0}^{n}
x_{\sigma^{-1}(1)}\t\cdots\t\underbrace{\{x_{\sigma^{-1}(j)},A\}}_{j\textrm{th place}}\t\cdots\t x_{\sigma^{-1}(n)}\\ &=& \sum_{i=0,
\sigma(j)=i}^{n} x_{\sigma^{-1}(1)}\t\cdots\t\underbrace{\{x_i,A\}}_{\sigma^{-1}(i)\textrm{th place}}\t\cdots\t x_{\sigma^{-1}(n)}\\
&=&\sigma\cdot\sum_{i=1}^{n}x_1\t\cdots\t \{x_i,A\}\t\cdots\t x_n\\ &=&\sigma\circ\rho(A)x \ .  \end{eqnarray*} And the proof is
completed.\hfill$\Box$\medskip

We then have to transform the module of chains $\Sp_{2m}(\Com)^{\t n}$ so that the $(V^{\t n})_{\Sp(\KK)}$ appears.
\begin{lemma}\label{prop:decomposition} Let $V_m$ be the symplectic $2m$-dimensional vector space.  The following equality holds :  $$
(\Sp_{2m}(\Com)^{\t n})_{\Sp_{2m}(\KK)}=\bigoplus_{\begin{array}{c} \scriptstyle k_1+\cdots+k_n=2r \\ \scriptstyle k_i\geq 2
\end{array}}(( \Vtm {2r})_{\Sp_{2m}(\KK)})_{\Sigma_{k_1}\times\cdots\times\Sigma_{k_n}} \ .  $$ \end{lemma}

\begin{proof} Direct calculation leads to the following equalities :  \begin{eqnarray*} &&\Sp_{2m}(\Com)^{\t n}:=S^{\geq 2}(V_m)^{\t n}=
(\bigoplus_{k\geq 2} (\Vtm k)_{\Sigma_k})^{\t n}\\ &&\qquad= \bigoplus_{\begin{array}{c} \scriptstyle (k_1,\cdots,k_n) \\ \scriptstyle
k_i\geq 2 \end{array}}( \Vtm {k_1})_{\Sigma_{k_1}}\t\cdots\t(\Vtm {k_n})_{\Sigma_{k_n}}\\ &&\qquad=\bigoplus_{\begin{array}{c}
\scriptstyle k_1+\cdots+k_n=2r \\ \scriptstyle k_i\geq 2 \end{array}} ( \Vtm {2r})_{\Sigma_{k_1}\times\cdots\times\Sigma_{k_n}}\oplus
\bigoplus_{\begin{array}{c} \scriptstyle k_1+\cdots+k_n=2r+1 \\ \scriptstyle k_i\geq 2 \end{array}} ( \Vtm
{2r+1})_{\Sigma_{k_1}\times\cdots\times\Sigma_{k_n}} \end{eqnarray*} Then, applying proposition \ref{actioncommute} to
$(\Sp_{2m}(\Com)^{\t n})_{\Sp_{2m}(\KK)}$ and proposition \ref{odd} completes the proof.  \end{proof}

We are now ready to apply the co-invariant theory to obtain the following proposition:  \begin{proposition}\label{prop:sptoAn} The chain
module $(\Sp(\Com))^{\t n}$ is isomorphic to $$ (\Sp(\Com)^{\t n})_{\Sp(\KK)}=\bigoplus_{\begin{array}{c} \scriptstyle k_1+\cdots+k_n=2r
\\ \scriptstyle k_i\geq 2 \end{array}}( A_r)_{\Sigma_{k_1}\times\cdots\times\Sigma_{k_n}} \ .  $$\end{proposition} \begin{proof} It
suffices to apply lemma \ref{prop:decomposition} and the co-invariant theory theorems \ref{theo:FFT} and \ref{theo:SFT}.  \end{proof}

\begin{remark}\label{rem:split} The map $T(\Sp(\Com))\to T(\Sp(Com))_{\Sp(\KK)}\to \oplus
(A_r)_{\Sigma_{k_1}\times\cdots\times\Sigma_{k_n}}$ admits a splitting $S$ induced by the following construction :  to a monomial
$y_{i_1,i_2}\cdots y_{i_{2r-1},i_{2r}}$ we associate a monomial in indeterminates $\{p_1,q_1,\ldots\}$ in $\Vt{2r}$ such that at the
place $i_{2k-1}$ there is a $p_k$ and at the place $i_{2k}$ there is a $q_k$.  Moreover this map is a $\Sigma_n$-morphism.  \end{remark}
\begin{example} The image by $S$ of $y_{14}y_{27}y_{35}y_{86}$ is $p_1p_2p_3q_1q_3q_4q_2p_4$.  \end{example}

\section{Second step :  Chord diagrams chain complex} First, we show that the vector space $A_r$ coming from the co-invariant theory is
isomorphic to the vector space spanned by the base pointed chord diagrams.  So we can understand the Leibniz chain complex coming from
the co-invariant theory as a chord diagram chain complex.

\begin{definition} Let $m$ be a positive integer.  A partition $c$ of $\{1,\ldots ,2m\}$ such that every $x \in c$ is a set consisting
precisely of two elements will be called a \emph{(base pointed) chord diagram}, where the base point is the set with the element $1$.
The set of all such chord diagrams will be denoted by $\Diag_m$.  A chord diagram with an ordering of each two point set will be called
an \emph{oriented chord diagram}.  The set of all oriented chord diagrams will be denoted by $\stackrel{\to}{\Diag}$.    There is a symmetric action on the chord diagram. Given an oriented chord diagram $c:=\{[i_1,j_1],\ldots,[i_m,j_m]\}$ and a permutation $\sigma
\in S_{2m}$, the symmetric group
$S_{2m}$ acts on $\stackrel{\to}{\Diag_m}$ as follows:
\[ \sigma\cdot c := \{[\sigma^{-1}(i_1),\sigma^{-1}(j_1)],\ldots,[\sigma^{-1}(i_m),\sigma^{-1}(j_m)]\}.  \] 

The symmetric group $S_{2m}$
acts on $\Diag_m$ in a similar way.  A \emph{labelled} diagram is a diagram $D\in\Diag_m$ together with a map between $\{1,\ldots , 2m
\}$ to the set of labellings $\{p_1,q_1,\ldots\}$.
 Any element in $\{1,\ldots 2m\}$ is called a \emph{vertex} and any set of two elements is called a \emph{label}.
\end{definition}

\begin{example}\label{ex:diagramgeometric}
There is a geometrical interpretation of a chord diagram :  it's
usual to put vertices on a circle and to draw each chord inside the disk.  For example the diagram
$D=\{\{1,4\},\{2,7\},\{3,5\},\{8,6\}\}\in \Diag_4$ can be geometrically represented as follows:  cf.\ figure \ref{figure:diagexample}.

\begin{figure}[htbp] 
\begin{center} 
\input{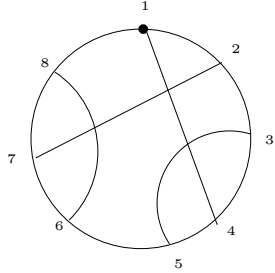} 
\caption{Geometric representation of the chord diagram D.}
\label{figure:diagexample} 
\end{center} 
\end{figure}

If the diagram is oriented then the orientation will be represented as some arrows
on the chords as in figure \ref{figure:orienteddiagexample}.  

\begin{figure}[htbp] 
\begin{center} 
\input{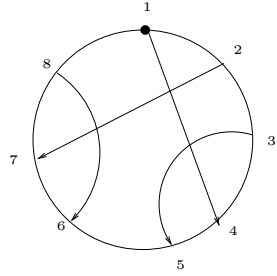}
\caption{Geometric representation of the oriented chord diagram D.}  
\label{figure:orienteddiagexample} 
\end{center}
\end{figure} 

The numbering on the circle is defined as follows :  determine a base point  and label it with  $1$ and label increasingly according to a positive orientation.
\end{example}

\begin{proposition}\label{prop:AnDiag}
 The space $A_m$ is isomorphic as a $\Sigma_n$-module to the vector space spanned by the collection
of chord diagrams with $m$ chords, denoted $\KK[\Diag_m]$.  
\end{proposition}

\begin{proof} 
Let $\beta:=y_{j_1,j_2}\ldots y_{j_{2m-1}j_{2m}}\in A_m$ be a monomial.  As any monomial in $A_m$ satisfies the relation
$y_{ij}=-y_{ji}$, $\beta$ can be rewritten (up to a sign) as the standard monomial $y_{i_1 i_2}\ldots y_{i_{2m-1} i_{2n}}$ with $i_1=1$,
$i_{2k-1}<i_{2k}$ and $i_{2k-1}<i_{2k+1}$ for all $k\in\{1,\ldots,m\}$.  (This remark will mod out the orientation on the diagrams.)
Therefore, we can focus on ordered monomials and construct the isomorphism.

Consider $\gamma:=y_{i_1 i_2}\ldots y_{i_{2m-1} i_{2m}}$ an ordered monomial.  It can be associated to the following base-pointed diagram
$\{\{1,i_2\},\ldots,\{\textrm{min}(i_{2m-1},i_{2m}),\textrm{max}(i_{2n-1},i_{2n})\}\}\in\Diag_m$.

By linearity, this construction gives rise to a map $\phi:A_m\to \KK[\Diag_m]$ which is a vector space isomorphism.  (The inverse map is
clear).

We verify that  the map $\phi$ is $\Sigma_{n}$-equivariant.  The permutation $\sigma\in \Sigma_{2n}$ acts on a standard monomial
as :$$ \sigma\cdot y_{i_1 i_2}\ldots y_{i_{2n-1} i_{2n}}= y_{\sigma^{-1}(i_1) \sigma^{-1}(i_2)}\ldots y_{\sigma^{-1}(i_{2n-1})
\sigma^{-1}(i_{2n})}\ .  $$ The chord diagram $\phi(\sigma\cdot y_{i_1 i_2}\ldots y_{i_{2n-1} i_{2n}})$ is defined as the following partition
:$$\{\{\sigma^{-1}(i_1) ,\sigma^{-1}(i_2)\},\ldots, \{\sigma^{-1}(i_{2n-1}), \sigma^{-1}(i_{2n})\}\}\ .$$ This is exactly
$\sigma\cdot\phi( y_{i_1 i_2}\ldots y_{i_{2n-1} i_{2n}})$.  Therefore $\phi$ is a $\Sigma_n$-module isomorphism.  
\end{proof}

\begin{example} 
The image of the monomial $y_{14}y_{27}y_{35}y_{86}$ under  the map $\phi$  is the diagram defined in example
\ref{ex:diagramgeometric}, which has the geometric representation of figure \ref{figure:diagexample}.  
\end{example}

Though it is unnecessary for the proof of theorem \ref{thm:LeibnizKont}, we will define the chain complex of chord diagrams and then
prove that this chain complex is quasi isomorphic to the chain complex defined for $\Sp(\Com)$.  This part can be skipped as in the next
section, we show that the chain complex $\Sp(\Com)$ is quasi-isomorphic to the chain complex defined on the graphs.

\begin{remark}\label{rem:package}
Let $\Gamma\in(\Diag)_{\Sigma_{k_{1}}\times\ldots\times\Sigma_{k_{n}}}$ be a diagram.We will call a \emph{package} the subset of vertices on which a $\Sigma_{k_{i}}$ acts. 
Let $D\in\Diagto$ be a diagram obtained by $D'\in\Diagto$ by a change of orientation on a chord $e=[i,j]$ into $[j,i]$. We set $D\sim -D'$.  

The algebra spanned by the oriented chord diagrams quotiented by the above equivalence relation 
is exactly the algebra spanned by the chord diagrams without any chord with two incident half edges in the same package $\Sigma_i$.  \end{remark}

\begin{definition} Let $D\in\Diag_m$ be a diagram and $e$ a chord of $D$.  We define the \emph{contraction} of $D$ by $e$ as a new
diagram denoted $D/e\in\Diag_{m-1}$ where the set $\{1,\ldots,2(m-1)\}$ admits for partition the shifting of the partition $c$ of $D$
where $e$ has been deleted (standardisation).  \end{definition}

\begin{example} Let $D$ be the diagram defined in example \ref{ex:diagramgeometric}, and let $e$ be the edge $\{2,7\}$.  The contraction
of $D$ by the edge $e$ is the following diagram $D/e= \{\{1,3\},\{2,4\},\{5,6\}\}$ which can be geometrically represented as figure
\ref{figure:diag_contraction}.  \begin{figure}[htbp] \begin{center} \input{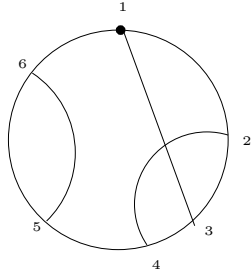} \caption{Geometric representation
of the chord diagram D/e.}  \label{figure:diag_contraction} \end{center} \end{figure} \end{example}

\begin{definition} Let $D\in \stackrel{\to}{\Diag}$ be an oriented chord diagram and let $c=[i,j]$ be an oriented edge of $D$.  We define the
sign $\epsilon(c)\in\{-1,1\}$ as follows:  \begin{equation*} \epsilon(c)=\left\{\begin{array}{ccc} -1&\textrm{if } i>j\ ,\\ 1&\textrm{if
}i<j\ .  \end{array}\right .  \end{equation*} \end{definition}

\begin{definition}

Denote by $C^m(\Diagto,\delta)$ the chain complex defined as follows.  The $n$th chain module is
$\oplus_{k_1+\cdots+k_r=n}(\Diagto_n)_{\Sigma_{k_1}\times\cdots\times\Sigma_{k_r}}$.  These chain modules are endowed with the following
differential :  let $D\in\Diagto_n$, 
we define:
$$\partial( \stackrel{\to} {D}):=\sum_{c}(-1)^i\epsilon(c) \stackrel{\to} {D/c}\ ,$$ 
where the sum runs over all the chords $c=[a,b]$ and $i$ is defined such as $\sum_{p=1}^{i}k_p\leq \textrm{max}(a,b)< \sum_{p=1}^{i+1}k_p$.
 
Remark that the result of this differential lives in
$$\oplus(\Diagto_{n-1})_{\Sigma_{k_1}\times\cdots\Sigma_{k_{i}+k_j-2}\times\cdots\hat{\Sigma_{k_j}}\times\cdots\times\Sigma_{k_r}}\ ,$$
where the sum is extended to all $[i,j]$ which are chords of the diagram.  
\end{definition} 

\begin{proposition}
The differential passes through the equivalence relation of remark \ref{rem:package}.
\end{proposition}

\begin{proof} The differential does not
depend on the representative oriented diagram of $D$.  Indeed, let $ \stackrel{\to} {D}$ be an oriented representative of the chord
diagram $D$, and let $c=[a,b]$ be one of its oriented chord.  Consider the oriented diagram $ \stackrel{\to} {D'}$ with the same
orientations as $ \stackrel{\to} {D}$ except for the chord linking $a$ to $b$ where we consider the orientation $[b,a]=:c'$.  A direct
computation ends the proof:  \begin{eqnarray*} \partial( \stackrel{\to} {D})-\partial( \stackrel{\to} {D'})&=&(-1)^i(\epsilon(c)[
\stackrel{\to} {D/c}]-\epsilon(c')[ \stackrel{\to} {D/c}])\\ &=&(-1)^i\epsilon(c)(1-1)[ \stackrel{\to} {D/c}]\\ &=&0\ .  \end{eqnarray*}
Indeed, any other changes in the orientation will just lead to a sum of the above equality.  
\end{proof}

\begin{example}\label{ex:diagdiff}Consider the chord diagram $D\in(\Diag)_{\Sigma_3\times\Sigma_3\times \Sigma_2}$ defined in example
\ref{ex:diagramgeometric}.  In order to compute its differential we will consider the representative oriented chord diagram $\hat
D\in(\stackrel{\to}{\Diag})_{\Sigma_3\times\Sigma_3\times \Sigma_2}$ also defined in the same example:  \begin{eqnarray*}
\partial(D)=\partial(\hat D)&=&(-1)^3\{\{1,3\},\{2,4\},\{6,5\}\}-(-1)^3\{\{1,4\},\{2,6\},\{3,5\}\}\ .  \end{eqnarray*} The differential
is more easily expressed thanks to a geometrical representation where we explicitly materialise the packages $\Sigma_{k_i}$ in the result
of the differential.  In figure \ref{figure:diag_differential}, we give the result of all contracted diagrams without taking into account
that the diagrams such that a chord is included in a package $\Sigma_{k_i}$ are null.  Then, we take this relation into account to
give the result of the differential in figure \ref{figure:diagexdiffresult}.

\begin{figure}[htbp] \begin{center} \input{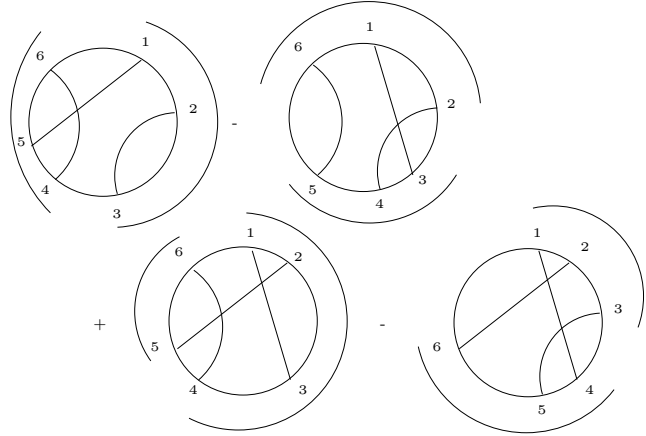} \caption{Geometric representation of the $\partial(D)$ without the
equivalence relation.}  \label{figure:diag_differential} \end{center} \end{figure}

\begin{figure}[htbp] \begin{center} \input{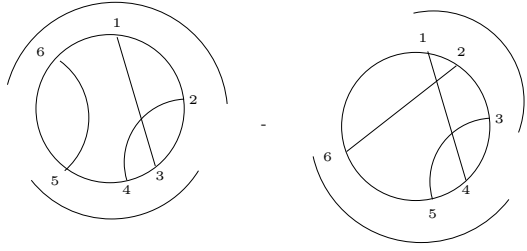} \caption{Geometric representation of the $\partial(D)$.}
\label{figure:diagexdiffresult} \end{center} \end{figure} \end{example} \begin{proposition} The differential on the chain complex of
diagrams is induced by the Leibniz differential of the chain complex of $\Sp(\Com)$.  \end{proposition} \begin{proof} See proof of
proposition \ref{prop:iso_homology}.  \end{proof}

\section{Second step:  Kontsevich idea and symmetric graph complex} In the context of Lie homology, Kontsevich's major contribution is
to give an isomorphism between the quotient space of chord diagrams by the action of the symmetric group
$(\KK[\Diag])_{\Sigma_{k_1}\times\cdots\times\Sigma_{k_n}}$ and some graphs.  We will mimic this idea to show that the graphs arising
are a symmetric version of Kontsevich graphs:  they admit an labelling of the vertices.  To avoid the problem of signs in the
differential, we will use a trick seen in \cite{CV}, and consider the chain complex of oriented symmetric graphs.

\begin{proposition}\label{prop:DiagGraph} Let $r,k_1,\ldots,k_n$ be integers satisfying the relation $\sum_{i=1}^n k_i=2r$ and such that
$k_i\geq 2$ for all $1\leq i\leq n$.  There exists a vector space isomorphism : $$(\KK[\Diag_r])_{\Sigma_{k_1}\times\cdots\times\Sigma_{k_{n}}}\cong \KK[\GG_{k_1\cdots k_n}]\ .$$  \end{proposition}

\begin{proof} We construct explicitly the map $\varphi:(\KK[\Diag_r])_{\Sigma_{k_1}\times\cdots\times\Sigma_{k_{n}}}\rightarrow \KK[\GG_{k_1\cdots k_n}]$. We associate to any chord diagram $D:=\{\{1,i_2\},\ldots,\{i_{2r-1},i_{2r}\}\}$  the following graph :  
\begin{eqnarray*}
V(G)&=&\{1,,\ldots,n\}\ ,\\ 
E(G)&=&\{\{1,\norm(i_2)\},\ldots,\{\norm(i_{2r-1}),\norm(i_{2r})\}\}\ ,
\end{eqnarray*}
where the map $\norm:\{1,\ldots, 2r\}\to\{1,\ldots,n\}$ is defined as follows:
\begin{eqnarray*} 
\norm(j)= l, \textrm{ where $l$ is defined by }
\sum_{i=1}^{l}k_i\leq j < \sum_{i=1}^{l+1}k_i 
\end{eqnarray*} 
We extend linearly this construction to define the map
$\varphi:(\KK[\Diag_r])_{\Sigma_{k_1}\times\cdots\times\Sigma_{k_{n}}}\rightarrow \KK[\GG_{k_1\cdots k_n}]$ .

This map admits an inverse map defined by the following algorithm.  Let $r:=\sum_{i=1}^n k_i$, let $G\in\GG_{k_1,\ldots,k_n}$ be a graph
with edges $(i_k,i_l)$.  We construct a diagram $D\in(\Diag_r)_{\Sigma_{k_1}\times\cdots\times\Sigma_{k_n}}$.  The algorithm to define
the edges is the following.  Let $\mathit{comp=0}$, $\mathit{ind}=1$ and $D'=\{\alpha_G(e_1)\cdots,\alpha_G(e_r)\}$.  Go through each set
(of cardinality two) of $D'$ if the element $i_x$ is equal to $\mathit{ind}$ then we indent this element into $i_x+\mathit{comp}$.  Then
$\mathit{comp}$ takes the value $\mathit{comp}+1$, and we repeat the algorithm.  When, the algorithm is over all the edges,
$\mathit{ind}$ takes the value $\mathit{ind}+1$.  (The counter $\mathit{comp}$ should be taking the value
$1+\sum_{i=1}^{\mathit{ind-1}}k_i$.) Restart the algorithm on the vertices that were not modified.

It is clear that the two maps are inverse to each other.  Therefore $\phi$ is an isomorphism, and thus ends the proof.  \end{proof}

\begin{remark} The fact that the diagrams with a chord included in  a package $\Sigma_i$ are null induces that the graphs with loops are null.
Indeed, let $G\in\GG_{k_1,\cdots,k_m}$ be a graph and $e=(i,i)$ be a loop of $G$.  Then $\phi\mn\circ\varphi\mn(G)$ is the monomial in
indeterminates $y_{kl}$ of the form $\cdots y_{j,j+1}\cdots\in(A_{r})_{\Sigma_{k_1}\times\cdots\times\Sigma_{k_m}}$.  A representative of this
monomial in $(\Vt{2r})_{\Sigma_{k_1}\times\cdots\times\Sigma_{k_m}}$ is a monomial such that $p_1$ is at the place $j$ and $q_1$ is at the
place $j+1$.  By the symmetric action, this monomial is also equal to the same monomial where $p_1$ is at the place $j+1$ and $q_1$ is at
the place $j$.  Therefore, going through $T^*$, it gives the monomial $\cdots y_{j+1,j}\cdots\in(\
A_{r})_{\Sigma_{k_1}\times\cdots\times\Sigma_{k_m}}$.  By the equivalence relation, this monomial is exactly $-\cdots y_{j,j+1}\cdots$.
Therefore, we proved that graphs with loops are null.  \end{remark}

\begin{example} The image of the diagram $D$ of the example \ref{ex:diagdiff} by $\phi$ is the following graph:
$$(\{1,2,3\},\{\{1,2\},\{1,3\},\{1,2\},\{2,3\}\}) $$ It can be understood geometrically by figure \ref{figure:diagramtograph}.
\begin{figure}\label{diagtograph} \begin{center} \input{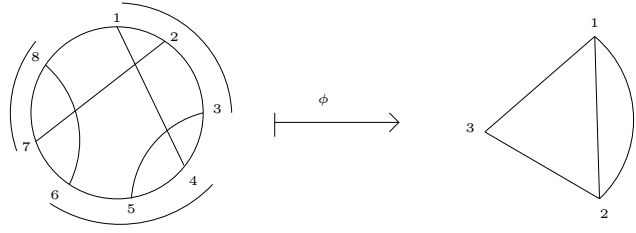} \caption{From chord diagram to graphs}
\label{figure:diagramtograph} \end{center} \end{figure} \end{example}

\begin{proposition}\label{prop:iso_homology} The Leibniz homology of the Lie algebra $\Sp(\Com)$ is isomorphic to the homology of the 
symmetric graph complex:  $$ HL_n(\Sp(\Com))=H_n(\KK[\GG],\delta) $$

\end{proposition}

\begin{proof} The first step of the proof of theorem \ref{thm:LeibnizKont} gives a quasi-isomorphism between the Leibniz chain complex of
$\Sp(\Com)$ and the Leibniz chain complex defined for $\Sp(Com)_{\Sp(\KK)}$.  Then, propositions \ref{prop:sptoAn}, \ref{prop:AnDiag} and
\ref{prop:DiagGraph} give a vector space isomorphism between the Leibniz chain complex of $(\Sp(\Com))_{\Sp(\KK)}$ and the vector space
of graphs, namely $\varphi\circ\phi\circ T^*$.  Therefore, it suffices to show that the differential defined on the graphs is exactly
induced by the differential of $\mathrm{CL}(\Sp(\Com))$.  In order to do so, we will consider a graph, and through splittings consider a
representative of this graph in $T(\Sp(\Com))$.  Then, we will explicit the differential of this representative, and we will give this
result in terms of graphs thanks to the vector space isomorphism.  Since the result is exactly the differential defined on the complex of
graphs, it will end the proof.

Consider the integers $k_1,\ldots,k_n$ and define $r:=\sum_{i=1}^n k_i$.  Let $G=(V(G),E(G),\alpha_G)$ be a graph in
$\GG_{k_1,\cdots,kn}$.

First, we construct a representative of $G$ in $T(\Sp(Com))$ as $S\circ\phi\mn\circ\varphi\mn(G)$.  The map $\varphi\mn$ in the proof of
proposition \ref{prop:DiagGraph} gives the construction of a diagram $D\in (\Diag_r)_{\Sigma_{k_1}\times\cdots\Sigma_{k_n}}$.  Then, to
have a representative of $D$ in $(A_r)_{\Sigma_{k_1}\times\cdots\Sigma_{k_n}}$, we consider the monomial $m:=\phi\mn(D)\ ,$ see the proof
of proposition \ref{prop:AnDiag}.  Moreover, this monomial admits a representative in $T(\Sp(\Com))$ defined as $P:=S(m)$ (see remark
\ref{rem:split}).  This element $P$ is of the form $\pm\underbrace{p_1 p_2\ldots
p_{k_1}}_{k_1}\t\underbrace{\ldots}_{k_2}\t\ldots\t\underbrace{q_{i_1}\ldots q_{i_{k_n}}}_{k_n}$.

Then, we explicit the differential of the representative.  Taking the differential of this element where each $p_i$ and each $q_i$ appear
once in this order in different factors of the tensor product gives the sum of signed elements of graduation $n-1$ where two factors have
merged and one couple in this merged factor is omitted.  Indeed, let us denote \begin{equation*} \tilde p_i:= \left\{ \begin{array}{rl}
p_i & \mbox{if } i \leq r \\ q_{i-r} & \mbox{si } i \geq r+1 \end{array} \right.  \ .  \end{equation*} Then, \begin{eqnarray*}
&&d(\underbrace{\tilde p_{i_1}\ldots\tilde p_{i_{k_1}}}_{k_1}\t\ldots\t \underbrace{\ldots p_{i_{2r}}}_{k_n})=\\
&&\sum^{n}_{\begin{array}{c} \scriptstyle j<k \\ \scriptstyle j=1,k=2 \end{array}}(-1)^{j}\underbrace{\tilde p_{i_1}\ldots\tilde
p_{i_{k_1}}}_{k_1}\t \ldots\t\{\underbrace{\tilde p_{i_l}\ldots \tilde p_{i_{l+k_i}}}_{k_i},\underbrace{\tilde p_{i_m}\ldots \tilde
p_{i_{m+k_j}}}_{k_j}\}\t \ldots\t \underbrace{\ldots p_{i_{2r}}}_{k_n}\ , \end{eqnarray*} where $l=\sum_{s=0}^{i-1}k_s$ and
$m=\sum_{s=0}^{j-1}k_s$.  The only way the element of the sum is non-trivial is that there exists at least a couple $p_s$ in the $k_i$th
factor of the tensor product and $q_s$ the $k_j$th factor (the case $p_s$ in the $k_j$th factor of the tensor product and $q_s$ in the
$k_i$th factor does not happen in our construction therefore no sign will appear from here).  As a couple $(p_s,q_s)$ appears only once it
is clear that the $k_i$ and the $k_j$ factor will concatenate omitting the couple $(p_s,q_s)$.  The sign that appears depends on the
number $j$ of the factor of $\Sp(\Com)^{\t n}$ where the element $q_i$ of the couple appears (as the couple appears in this order in the
tensor factors).

The result of the differential can be understood in terms of graphs thanks to the vector space isomorphism $\varphi\circ\phi\circ
T^\star$.  The number $j$ is represented as the $j$th vertex of the graph.  Moreover, the omission of the couple is translated by the
disappearance of the appropriate vertex and the identification of its two edges giving rise to the necessary changes of vertices which is
exactly the contraction of the graph with this vertex.  And so, we can conclude that :  $$\delta(G)=\sum_{e=[i,j]\in
E(G)}(-1)^{\textrm{max}(i,j)} \epsilon(i,j) G/e \textrm{ for all }G\in\GG\ .$$ (In our construction, the oriented representative that we
take is exactly the one such that every $\epsilon(i,j)=1$, and therefore $\max(i,j)=j$.  )

It is exactly the differential defined on the chain complex of graphs.  This ends the proof.  \end{proof}

\begin{example} We sketch the idea of the proof on an example.  Let $G$ be the graph defined as
$(\{1,2,3\},\{\{1,2\},\{1,3\},\{1,2\},\{2,3\}\})$ with geometric interpretation as in figure \ref{figure:graphexample}.  This graph can
be lifted into the following diagram $\{\{1,4\},\{2,5\},\{3,7\},\{6,8\}\}\in(\Diag_4)_{\Sigma_3\times\Sigma_3\times\Sigma_2}$ which is
geometrically interperted as the diagram in figure \ref{figure:diagexample}.  This diagram  is isomorphic to the monomial
$y_{14}y_{27}y_{35}y_{68}$ in $(A_4)_{\Sigma_3\times\Sigma_3\times\Sigma_2}$.  It can be lifted as the following monomial in
$T(\Sp(\Com))$ :  $$P:=p_1p_2p_3\t q_1q_2p_4\t q_3q_4 $$ Taking the differential of this element gives the following result:
\begin{eqnarray*} d(p_1p_2p_3\t q_1q_2p_4\t q_3q_4)&=&p_2p_3q_2p_4\t q_3q_4+p_1p_3 q_1p_4\t q_3q_4\\&& -p_1p_2p_3\t q_1q_2q_3-p_1p_2q_4\t
q_1q_2p_4\ .  \end{eqnarray*} To interpret this result in terms of graphs, we compute $\varphi\circ\phi\circ T^*$ of the result.  By the
map $T^*$ we obtain the sum of the following monomials in $(A_3)_{\Sigma_{4}\times\Sigma_{2}}\oplus(A_3)_{\Sigma_{3}\times\Sigma_{3}}$ :
\begin{eqnarray*} &&2y_{13}y_{25}y_{46}-2y_{14}y_{25}y_{36}\ ,\\ &=& -2y_{14}y_{25}y_{36} \end{eqnarray*} by the symmetric action, then
the map $\phi$ gives the following diagram in $(\Diag_3)_{\Sigma_{3}\times\Sigma_{3}}$ :  $$(-2(\{1,4\},\{2,5\},\{3,6\}))\ , $$ see
figure \ref{figure:diagobtdiff}, finally the map $\varphi$ gives the sum of graphs :  $$-2(\{1,2\},\{\{1,2\},\{1,2\},\{1,2\}\}) $$ see
figure \ref{figure:graphdiffresult}.  By example \ref{ex:graphdiff} we realise that the two calculations of the differential are
identical.

\begin{figure}[htbp] \begin{center} \input{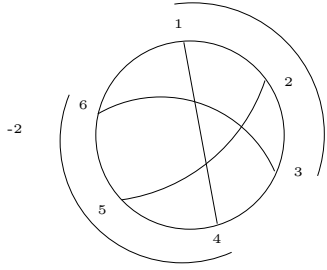} \caption{Geometrical interpretation of the diagram $\phi\circ T^*(P)$}
\label{figure:diagobtdiff}\end{center} \end{figure} \end{example}

\begin{proposition}\label{prop:isocon}
There exists a vector space isomorphism between the Leibniz homology of $\Sp(\Com)$ and the chain complex of connected graphs  :
$$HL_{n}(\Sp(\Com))\cong T(H_{n}(\KK[\GG_{c}])) \ .$$
\end{proposition}

\begin{proof} 
The above proposition ensures that $HL_n(\Sp(\Com))\cong H_n(\KK[\GG],\delta)$.
Moreover there is a vector space isomorphism between the tensor module over the vector space spanned by connected graphs and the vector space spanned by graphs :
$HL_{n}(\Sp(\Com))\cong H_{n}(T(\KK[\GG_{c}]))$. It is well-known, cf.\ appendix B of Quillen \cite{Qa}, that the functor $T$ and $H$ commute. This ends the proof.
\end{proof}

\section{Third step: explicit homotopies}
This step reduces the computation of the homology of the vector space spanned by connected graphs. As in the
Kontsevich case we can reduce the computation of the homology to the complexes of graphs
which have no bivalent vertex. To avoid a spectral sequence, we show the acyclicity of
the quotient complex of the connected graphs by the polygons and the graphs with no
bivalent vertex, by producing an explicit  homotopy.
Moreover, with few changes, this homotopy could be used in the Lie context.

Let $\KK[\GG_{\rm{c}}]$ be the subcomplex of connected graphs. We denote by
$\KK[\GG_{\rm{c}}^3]$ the subcomplex of graphs with no bivalent vertex. The
subcomplex of graphs with only bivalent vertices are the polygons and is denoted by $\KK[P]$.

\begin{proposition}\label{prop:acyclic}
The subcomplex of graphs $\KK[\GG_{\rm{c}}/_{P\oplus\GG_c^3}]$ is acyclic.
\end{proposition} 

\begin{proof}
We construct a homotopy $h:\KK[\GG_{\rm{c}}/_{P\oplus\GG_{c}^3}]\to
\KK[\GG_{c}/_{P\oplus\GG_{c}^3}]$.
Let $L_k$ be the ladder graph with $k$ bivalent vertices. Let $G$ be a connected graph
with $n$ vertices and with $m$ ladders. We define $G_{+i}$ to be the graph $G$ where the
ladder $i$ with $k$ inner vertices  $L_k$ is replaced by $L_{k+1}$ such that the added
vertex is the last one and that it is labelled with $n+1$.

$$h(G):=\sum_{i}\frac{(-1)^{n+1}}{m}G_{+i}\ .$$

We  verify that $hd+dh=\Id$. There are two cases to go through. The first one is when $h$
and $d$ are adding and contracting  edges of the same ladder:
see figure \ref{figure:preuve}.
\begin{figure}[htbp]
\begin{center}
\input{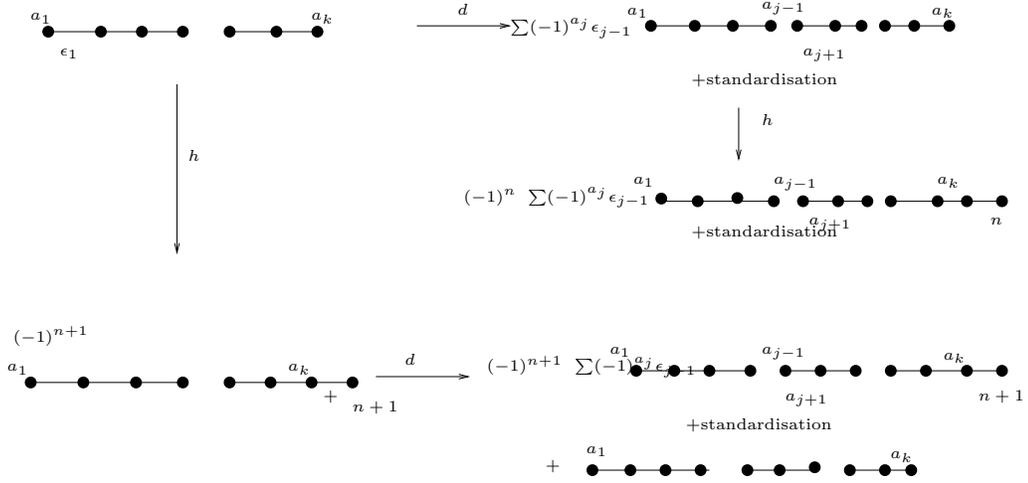}
\caption{Proof of the homotopy from $\Id$ to $0$.}
\label{figure:preuve}
\end{center}
\end{figure}

And therefore $dh+hd=\Id$.
The second is when  $d$ contractes an edge which is not part of this ladder. The two
actions anti-commute (because the sign of $h$ depends on the number of vertices which
fall by one with the differential) and therefore $dh+hd=0$ in this case. \end{proof}

\begin{proposition}\label{prop:polygon}
The homology of the subcomplex of graphs with only bivalent vertices is acyclic:
$$HL_*(\KK[P])\cong 0 \ .$$
\end{proposition}
\begin{proof}
The complex of the symplectic Lie algebra $T(S^2(V))$ is quasi-isomorphic to the subcomplex
of polygons. By proposition \ref{prop:sp2}, $S^2(V)$ is isomorphic to $\Sp(\KK)$. By
Pirashvili \cite{Pi}, the Leibniz homology of $\Sp(\KK)$ is null. This ends the proof.
\end{proof}

\begin{proposition}\label{prop:isovectorspace}
There exists a vector space isomorphism such that :
$$HL_{n}(\Sp(\Com)\cong T(H(\KK[\GG_{c}^3])) \ .$$
\end{proposition}

\begin{proof}
By proposition \ref{prop:isocon}, we need to show that $H(\KK[\GG_{c}])$ is isomorphic to $H(\KK[\GG_{c}^3])$. Moreover $\KK[\GG_{c}]$ is the sum of the subcomplexes $\KK[\GG_c^3]\oplus \KK[P] \oplus\KK[\GG_{c}/_{P\oplus\GG^3_c}]$. So, the above propositions end the proof.
\end{proof}

\section{Fourth step :  Graded differential Zinbiel-associative bialgebra}

We would like to prove that the isomorphism in proposition \ref{prop:isovectorspace} is not only a vector space isomorphism  but a Zinbiel-associative bialgebra isomorphism.
In order to do so, we must endow both homologies with this bialgebra structure. 

The Leibniz homology of any Leibniz algebra admits naturally a Zinbiel (the Leibniz Koszul dual) coalgebra structure. The associative operation is  particular to $HL_{*}(\Sp(\Com))$, and it is induced by the sum of the symplectic matrices.  Then, we show that the Zinbiel-associative structure on graphs, defined in definitions \ref{def:operationgraph} and \ref{def:cooperationgraph}, is induced by the Zinbiel-associative structure on the Leibniz homology $HL_{*}(\Sp(\Com))$, giving rise to the Zinbiel-associative isomorphism:
$$ HL_{*}(\Sp(\Com))\cong H_{*}(\KK[\GG]) \ .$$
Moreover, we can state a rigidity theorem, analogous to the Hopf-Borel theorem for co-commutative and commutative bialgebras, stating that a connected Zinbiel-associative bialgebra can be reconstructed from its primitives (see appendix). Therefore, we get the Zinbiel-associative isomorphism:
$$ H_{*}(\GG)\cong T(H_{*}(\GG_{c}))  \ .$$
Then, the last step is clear as the subcomplexes we consider ($P$, $\GG_{c}^3$ and $\GG_{c}/_{P\oplus\GG^3_c})$ are Zinbiel-associative subcomplexes. So we have a Zinbiel-associative isomorphism:
$$ HL_{*}(\Sp(\Com))\cong H_{*}(\KK[\GG_{c}^3]) \ .$$

\subsection{The Zinbiel coalgebra structure on a Leibniz chain complex $\mathrm{CL}_*(\gg)$}\label{sec:Zinbcooperation}

Let $\gg$ be a Leibniz algebra.  In his thesis \cite{O}, J.-M.  Oudom showed that the diagonal map $\gg\to\gg\times\gg:x\mapsto(x,x)$
induces a coproduct in the Leibniz chain complex of $\gg$, notably $(T\gg,\partial)$.  Indeed, the diagonal map induces a differential
map $\Phi:(T\gg,\partial)\to(T(\gg\times\gg),\partial)$.  The Zinbiel coproduct is defined as the projection of $\Phi$ on the first
component $T\gg\t T\gg\subset T(\gg\times\gg)$.  Moreover, J.-M.  Oudom showed that this differential map is exactly the co-half shuffle.
\begin{proposition}\cite{O} Let $\gg$ be a Leibniz algebra.  The co-half shuffle endows the Leibniz chain complex $T\gg$ with a
differential Zinbiel coalgebra structure with:  $$\Delta(g_1\ldots g_n):=g_1\sum_{p+q=n}\sum_{\underline i \in Sh_{p,q}}g_{i_1}\ldots
g_{i_p}\t g_{i_{p+1}}\ldots g_{i_n}\ ,$$ where the sum is extended over all $(p,q)$-shuffles $\underline i$ (i.e.  the in multi-indices
$\underline i=(i_1,\ldots,i_n)$ the integers $1, \ldots, p$ are ordered such as $p+1, \ldots, n$).  Moreover this Zinbiel coalgebra
structure is preserved on the Leibniz homology.\flushright{$\Box$} \end{proposition}

\subsection{The associative algebra structure the Leibniz chain complex $\mathrm{CL}_*(\Sp(\Com)_{\Sp(\KK)})$}\label{sec:assoperation}

To define the multiplication map, we consider the sum of matrices defined as :  $$\oplus :  \Sp(\Com)\times \Sp(\Com)\to
\Sp(\Com):(x,y)\mapsto E(x)+O(y) \ ,$$ where the maps $E:\Sp(\Com)\to \Sp(\Com)$ and $O:\Sp(\Com)\to \Sp(\Com)$ are induced by:
 \begin{eqnarray*} 
 E(p_i):=p_{2i},\ E(q_i):=q_{2i}\ ,\\
  O(p_i):=p_{2i-1},\ O(q_i):=q_{2i-1}\ ,
 \end{eqnarray*} see for example \cite{CV}.
 This maps induces
an operation on the chain Leibniz complex, by considering the injection of the first component $T(\Sp(\Com))\t T(\Sp(\Com))\subset
T(\Sp(\Com)\times\Sp(\Com))$.  It can moreover be shown that this map is associative on the Leibniz chain complex of $\Sp(\Com)_{\Sp(\KK)}$, see the proof of proposition
\ref{prop:associative-alg}.

Moreover, the Leibniz chain complex $\mathrm{CL}_*(\Sp(\Com)_{\Sp(\KK)})$ and the Leibniz homology $\mathrm{HL}_*(\Sp(\Com)_{\Sp(\KK)})$
admit a structure of Zinbiel-associative bialgebra.  It is proven thanks to the vector space isomorphism with the chain complex of
graphs.

\subsection{Zinbiel-associative bialgebra structure on the complex of graphs} The co-half shuffle endows $C^n(\GG,\delta)$ with a
structure of Zinbiel-associative bialgebra.

\begin{proposition}\label{prop:coprodinduct} The co-half shuffle defined on the graphs is induced by the co-half shuffle on the Leibniz
complex on $\Sp(\Com)$.  \end{proposition} \begin{proof}

The proof will be done in two steps. We will first focus on connected graphs. Let $G$ be
a connected graph. This graph can be seen as the inverse image of an element $w_1\cdots
w_n$ of $\Sp(\Com)^{\t n} $ as described in the proof of proposition 
\ref{prop:iso_homology} where $w_i\in V^{\t k_i}$ for certain $k_i$. First, we will apply
the co-half shuffle to this element and obtain an element of $\Sp(\Com)^{\t 2}$. Then, we
have to see this element as a graph once again under the map $(\phi\t\phi)\circ(
T^\star\t T^\star)$ which will  first  give rise to a diagram then, taking into account
the action of the cartesian product of symmetric groups leads to the graph. So, we have
to determine the non-zero elements rising from the map $T^\star\t T^\star$. It's
elementary to see that the  element $G\t 1$ will rise.  It is the only element. Indeed,
suppose  $T^\star (w_1w_{_2}\cdots w_{i_p})\neq 0$ and
$T^\star (w_{i_{p+1}}\cdots w_{i_n})\neq 0$. This induces that $w_1w_{i_2}\cdots
w_{i_p}\in V^{\t 2k_1}$ and  $w_{i_{p+1}}\cdots w_{i_n}\in V^{\t 2k_2}$ for $k_i\in \NN$.
And moreover,  supposes that there exists permutations $\sigma_1\in S_{2k_1}$ and
$\sigma_2\in S_{2k_2}$ such that $\omega^{\t k_1}(w_{\sigma(1)}w_{\sigma_1(i_{2})}\cdots
w_{\sigma_1({i_p})})=\pm 1$ and $\omega^{\t k_2}(w_{\sigma_2(i_{p+1})}\cdots
w_{\sigma_2({i_n})})=\pm 1$. Thus there exists a permutation 
$(\sigma_1(1),\cdots,\sigma_1(i_p),\sigma_2(i_{p+1}),\cdots,\sigma_2(i_n))\in S_n$ such
that $\omega^{\t
k_1+k_2}(\sigma_1(1),\cdots,\sigma_1(i_p),\sigma_2(i_{p+1}),\cdots,\sigma_2(i_n)) = \pm
1$. And by taking into account the symmetric group action this gives rise to a
non-connected graph. This implies that the connected and non-connected graphs are
isomorphic, which is a contradiction. Therefore $\Delta_<(G)=G\t 1$ for connected graphs.

The second part of the proof is done with the same arguments as above, considering a
non-connected graph. Indeed a non-connected graph is the disjoint union of connected
graphs.
\end{proof}
To illustrate the proof we consider the following example :
\begin{example}\label{ex:coproduitremonter}
The graph $G=(\{1,2,3\},\{\{1,2\},\{1,2\},\{1,3\},\{2,3\}\})$ considered in example
\ref{ex:graphexample}, can be lifted as $p_1p_2p_3\t q_2p_4\t q_1q_3q_4$ in
$T(\Sp(\Com))$. Taking it's co-half shuffle gives:
\begin{eqnarray*}
&&\DD(p_1p_2p_3\t q_2p_4\t q_1q_3q_4)=p_1p_2p_3\t q_2p_4\t q_1q_3q_4\bigotimes 1 +\\
&&p_1p_2p_3\t q_2p_4\bigotimes q_1q_3q_4+p_1p_2p_3\t q_1q_3q_4\bigotimes q_2p_4
+p_1p_2p_3\bigotimes q_2p_4\t q_1q_3q_4
\end{eqnarray*}
And applying $(\varphi\circ\phi\circ T^\star)\t (\varphi\circ\phi\circ T^\star)$ gives:
\begin{eqnarray*}
G\t 1
\end{eqnarray*}
which is exactly $\DD(G)$.

The non-connected graph $H=H_1\cdot H_1$ considered in example \ref{ex:nonconnexample}
admits for representative the following polynomial in $T(\Sp(\Com))$:
$$P:=p_1p_2\t q_1q_2\t p_3p_4\t q_3q_4\ .$$
Taking the differential of this element gives the following :
\begin{eqnarray*}
&&\DD(p_1p_2\t q_1q_2\t p_3p_4\t q_3q_4)=p_1p_2\t q_1q_2\t p_3p_4\t q_3q_4\bigotimes
1\\&&\qquad
+p_1p_2\bigotimes q_1q_2\t p_3p_4\t q_3q_4
 +p_1p_2\t q_1q_2\bigotimes p_3p_4\t q_3q_4\\&&\qquad+p_1p_2\t p_3p_4\bigotimes  q_1q_2\t
q_3q_4+
p_1p_2\t q_3q_4\bigotimes q_1q_2\t p_3p_4\\&&\qquad +p_1p_2\t q_1q_2\t p_3p_4\bigotimes
q_3q_4
+p_1p_2\t q_1q_2\t q_3q_4\bigotimes p_3p_4\\&&\qquad+
p_1p_2\t  p_3p_4\t q_3q_4\bigotimes q_1q_2 \ .
\end{eqnarray*}
To have the result in terms of graphs, we apply $(\varphi\circ\phi\circ T^\star)\t
(\varphi\circ\phi\circ T^\star)$ to obtain :
$$H\t 1+ 2H_1\t H_1 $$
which is exactly $\DD(H)$.
\end{example}
\begin{proposition}\label{prop:associative-alg}
The associative product, ordered disjoint union, on the complex of graphs is induced by
the associative structure on the Leibniz chain complex $(T(\Sp(\Com)))_\Sp$.
\end{proposition}
\begin{proof}
To prove this property, we will show that the ordered disjoint union is induced by the
operation on the Leibniz chain complex. The associativity on $\KK[\GG]$ is clear since
the operation is the ordered disjoint union of graphs. The associativity of the product
defined on $(T(\Sp(\Com)))_\Sp$ follows from the fact that  $T(\Sp(\Com))_\Sp$ is
isomorphic as vector space to $\KK[\GG]$.

We focus into proving that the ordered disjoint union of graphs is induced by
$T(\Sp(Com))\t T(\Sp(Com))\to T(\Sp(\Com))$. Let $G_1$ and $G_2$ be two graphs. These
graphs admit chord diagram representatives as constructed  in the proof of proposition
\ref{prop:DiagGraph}. Furthermore, by $\phi\mn$ these chord diagrams can be seen as  a sum
of monomials in variables $y_{ij}$. Last but not least, these monomials can be lifted up
as a polynomial $F_i$ in $T(\Sp(\Com))$, for $i=1,2$, by the split $S$ defined in remark
\ref{rem:split}.
These two polynomials can be seen as included in $T(\Sp(\Com)\times\Sp(\Com))$ by
decorating the variables of $F_1(p_1,q_1,\ldots)$ by $'$ and those of
$F_2(p_1,q_1,\ldots)$ by $''$. Then apply the operation $\oplus$ to them to obtain  the
following polynomial $F_1(p_2,q_2,\ldots,p_{2i},q_{2i},\ldots)\t
F_2(p_1,q_1,\ldots,p_{2i-1},q_{2i-1},\ldots)$. Then, by going through the isomorphism
$(\varphi\circ\phi\circ T^\star)\t (\varphi\circ\phi\circ  T^\star)$ we obtain the
ordered disjoint union of the graphs. Indeed, a vertex links variables $p$ and $q$ of
same indices, that is to say it links a $p_i$ with a $q_i$. Therefore the shifting we did
does not influence the vertices. Moreover it does not interfer in the decoration of the graph
as the second graph will be decorated with numbers following those from the first graph.
\end{proof}
To ease the comprehension of the proof, we consider the following example :
\begin{example}
Let $G$ and $H$ be the graphs of the above example \ref{ex:coproduitremonter}. The
associative product on these graphs is induced by the associative product on
$(T(\Sp(\Com)))_{\Sp(\KK)}$.  Indeed, by example \ref{ex:coproduitremonter} the two
graphs admit representatives in $(T(\Sp(\Com)))_{\Sp(\KK)}$. Apply the product
$\bigoplus$ to these representatives gives the following:
\begin{eqnarray*}
p_2p_4p_6\t q_4p_8\t q_2q_6q_8\t p_1p_3\t q_1q_3\t p_5p_7\t q_5q_7\ .
\end{eqnarray*}
By $(\varphi\circ\phi\circ T^\star)\t (\varphi\circ\phi\circ T^\star)$ we obtain the
result in terms of graphs:
$$(\{1,\ldots,7\},\{\{1,2\},\{1,2\},\{2,3\},\{1,3\},\{4,5\},\{4,5\},\{6,7\},\{6,7\}\})\
,$$
which is exactly $G\cdot H$.
\end{example}

\begin{proposition}\label{prop:Interchange_d_mu}
The induced product on the graph homology is associative and it is induced by the associative product on the Leibniz homology of $\Sp(\Com)_{\Sp(\KK)}$. 
\end{proposition}
\begin{proof}
First, we show that the operation on the homology of graphs is associative.
But, it is clear that $d\circ \mu-\mu\circ(\Id\t d+d\t\Id)=0$ on the graph complex, proving the associativity of the induced operation.
Then, we show that the induced operation on the Leibniz homology of $\Sp(\Com)_{\Sp(\KK)}$ is associative. Remark that on $T(\Sp(\Com))$, the following holds: 
\begin{eqnarray*}
d\circ \mu-\mu\circ(\Id\t d+d\t\Id)(v_{1}\ldots v_{p}\t v_{p+1}\ldots v_{p+q})=\\
 \sum_{1\leq i\leq p,p+1\leq j\leq p+q}v_{1}\ldots\{v_{i},v_{j}\}\ldots\hat{v_{j}}\ldots v_{p+q}\ \end{eqnarray*}
To prove the assertion, it suffices to remark that $T^*\circ(d\circ \mu-\mu\circ(\Id\t d+d\t\Id)))=0$ thanks to the vector space isomorphism with the vector space spanned by graphs.
\end{proof}

\begin{proposition}
There is a Zinbiel-associative bialgebra isomorphism between the Leibniz homology of
$\Sp(\Com)$ and the homology of graphs : $$HL_*(\Sp(\Com))\cong T(HL_*(\GG_c^3)) \ .$$
\end{proposition}

\begin{proof}
First we need to show that a Zinbiel-associative structure on the Leibniz chain complex
passes through homology. This is the case thanks to proposition \ref{prop:Interchange_d_D}
for the coproduct and proposition \ref{prop:Interchange_d_mu} for the product.

To define a Zinbiel-associative bialgebra structure on the Leibniz homology of
$\Sp(\Com)$ it suffices to consider the operation and co-operation induced through the
vector space isomorphism
 $$HL_*(\Sp(\Com))\cong HL_*(\Sp(\Com)_{\Sp(\KK)})\ ,$$
due to the Koszul trick.

Then,   propositions \ref{prop:coprodinduct} and \ref{prop:associative-alg} produce a
Zinbiel-associative bialgebra isomorphism :
$$CL_*(\Sp(\Com))\cong C_*(\KK[\GG])\ .$$

Then, apply the rigidity theorem \ref{theo:Zinb-as} to the connected  Zinbiel-associative bialgebra $C_*\KK[\GG]$:
$$   C_*(\KK[\GG])\cong T(\Prim \KK[\GG])\ .$$ 
By proposition \ref{prop:primgraph} the primitive graphs are the  connected graphs.
To conclude it suffices to realise that the subcomplexes considered in the third step are Zinbiel-associative subcomplexes. 
\end{proof}

\subsection{Proof of the Kontsevich theorem}
In this section, we give a short proof of Kontsevich's theorem in the flavour of the proof given in the Leibniz context. 

The set of graphs that Kontsevich considers is the set of classes of symmetric graphs quotiented by the signed symmetric action, that we denote $\GGG$.
We denote $\GGG_{c}$ the set of connected graphs, and $\GGG_{c}^3$ the set of connected graphs such that the vertices are of valency at least $3$. These graphs are geometrically the same as those in the Leibniz context, but without the numbering on the vertices.

Kontsevich theorem is stated as follows :
\begin{theorem}There exists a canonical co-commutative commutative bialgebra isomorphism :
$$H(\Sp(\Com)\cong \Lambda(H(\GGG_{c}^3)) \ .$$
\end{theorem}

The skeleton of the proof is as follows: 

First step, quotient the Chevalley-Eilenberg chain complex by the action of the reductive algebra $\Sp(\KK)$, thanks to the Koszul trick. Then apply the co-invariant theory to reduce the chain complex to the chain complex of chord diagrams (quotiented by the symmetric action). Then,  Kontsevich's idea is to consider the graphs, to code the quotient of the chord diagrams. The computation of the homology can be reduced to the computation of the homology of the connected graphs, which can be moreover reduced thanks to explicit homotopy.

The Chevalley-Eilenberg chain complex of $\Sp(\Com)$ is quasi-isomorphisc to the Chevalley-Eilenberg chain complex of $\Sp(\Com)_{\Sp(\KK)}$ similarly to proposition \ref{prop:quasiiso}.
$$ H(\Sp(\Com))\cong H(\Sp(\Com)_{\Sp(\KK)})\ . $$  
The co-invariant theory and direct computation gives a vector space isomorphism analogously to proposition \ref{prop:sptoAn}:
$$ (\Lambda^n(\Sp(\Com)))_{\Sp(\KK)}=\bigoplus_{\begin{array}{c} \scriptstyle k_1+\cdots+k_n=2r
\\ \scriptstyle k_i\geq 2 \end{array}}(( A_r)_{\Sigma_{k_1}\times\cdots\times\Sigma_{k_n}} )_{\Sigma_{n}}\ .  $$ 

Propositions \ref{prop:AnDiag} and \ref{prop:DiagGraph} still hold. Therefore, $A_{r}$ is isomorphic to the vector space spanned by chord diagrams, which quotiented by the symetric action is isomorphic to the vector space spanned by graphs $\GG$. It suffices to quotient the vector space of graphs $\GG$ by the symmetric action to conclude the existence of a vector space isomorphism:
\begin{equation}\label{eq:isograph}
H(\Sp(\Com))\cong H(\GGG) \ . 
\end{equation}
Any graph in $\GGG$ is a union of connected graphs, and the compatibility to the differential forces the existence of the following  vector space isomorphism :
   
\begin{equation}\label{eq:isograph2}
H(\Sp(\Com))\cong \Lambda(H(\GGG_{c})) \ . 
\end{equation}
Similarly to proposition \ref{prop:acyclic} and \ref{prop:polygon} the homology of the primitives can be reduced ti the vector space of connected graphs with no bivalent vertices $\GGG_{c}^3$.

This isomorphism is shown to be a co-commutative commutative bialgebra isomorphism as follows.
The chain complex of $\Sp(\Com)_{\Sp(\KK)}$ admits a  commutative and co-commutative bialgebra structure on the chain complex of graphs. The commutative operation on the Chevalley-Eilenberg complex of $\Sp(\Com)_{\Sp(\KK)}$ is induced by the sum of matrices, see section \ref{sec:assoperation}. The diagonal map induces the co-commutative co-operation.
The vector space isomorphism $C(\Sp(\Com))\cong C(\KK[\GGG])$ induces a structure of commutative co-commutative bialgebra structure on the chain complex of graphs. Therefore the isomorphism of equation \ref{eq:isograph} is a co-commutative commutatve bialgebra isomorphism. Moreover by the Hopf-Borel theorem, this connected commutative co-commutative bialgebra $\KK[\GGG]$ is isomorphic to the bialgebra $\Lambda(\Prim \GGG)$, where $\Prim \GGG= \GGG_{c}$.
Then, to conclude it suffices to realise that the subcomplexes considered in the last step, namely the subcomplex on polygons, the subcomplex on graphs with at least a bivalent vertex, and the subcomplex on graphs $\KK[\GGG_{c}^3]$ are bialgebra subcomplexes.
\begin{acknowledgement}
I would like to thank Jean-Louis Loday and Alain Brugui\`eres for their advisory. I am indebted to Belkacem Bendiffalah for introducing me to spectral sequences, and to Jean-Michel Oudom for interesting discussions and especially for the section ten.
This work has been partially supported by the ANR (Agence Nationnale de la Recherche) and by the project ECOS-Sud (Evaluation - Orientation de la Coopération Scientifique Sud).
\end{acknowledgement}

\newpage
\appendix
\section*{Appendix: Associative-Zinbiel bialgebras \cite{B}}
There is a celebrated theorem for classsical bialgebras known as the Milnor-Moore theorem
which states that a connected co-commutative bialgebra can be reconstructed thanks to its
primitive part. The goal of this appendix is to give an analogue of this theorem for
connected Zinbiel-associative bialgebras and dually for connected associative-Zinbiel
bialgebras.
\section{The Associative-Zinbiel structure theorem}\label{App:As-Zinb}
\subsection{Zinbiel algebra \cite{LZin}}
\begin{definition}
A Zinbiel algebra is a vector space $A$ endowed with a bilinear operation $\prec:A\t
A\to A$ verifying the following relation:
$$(x \prec y) \prec z = x \prec (y \prec z) + x\prec (z \prec y) \ , \ \  \forall x,y,z
\in A \ .$$
Moreover a Zinbiel algebra is said to be unital if it admits an element $1$ such that for
all $x \in A$ the following is verified:
\begin{equation}\label{unit}
\Bigg\{ \begin{array}{lcl}
1 \prec x &=& 0 \\
x \prec 1 &=& x \ , \ \ \forall x \in A \ ,
\end{array} \end{equation}
Note that $1\prec 1$ is not defined.
\end{definition}
Remark that the operation $\ast : A \times A \rightarrow A : (x,y) \longmapsto x\prec y +
y \prec x$
is associative, commutative and unital.

\begin{definition}
Let $A_{0}$ be a Zinbiel algebra. This algebra is  \emph{free over the vector space $V$},
if it satisfies the following universal property. Any map  $f : V\to A$, where $A$ is any
Zinbiel algebra,
extends uniquely into a Zinbiel morphism
$\tilde{f} : A_{0} \to A$. This can be summarised in the commutation of the following
diagram:
\begin{displaymath}
\xymatrix{
V \ar[r]^{i} \ar[dr]_{f}&A_{0} \ar[d]^{\tilde f}\\
& A&  \ .}
\end{displaymath}
\end{definition}
\begin{definition}
The shuffle algebra is the tensor module $T(V)$ over the vector space $V$ endowed with
the following operation $\Ss:T(V)\t T(V)\to T(V)$ defined as:
\begin{equation*}
v_1\cdots v_p \Ss v_{p+1} \cdots v_n := \sum_{\underline i \in Sh_{p,q}}v_{i_1}\ldots
v_{i_n} \in V^{\otimes n} \ ,
\end{equation*}
where the sum is extended to the $(p,q)$-shuffles $\underline i$, i.e. the multi-indice
$\underline i=(i_1,\ldots,i_n)$ has the property that  $1, \ldots, p$ are in this order
and so are $p+1, \ldots, n$. \end{definition}

\begin{proposition}\label{prop:Zlibre}
The free Zinbiel algebra over the vector space $V$, denoted $Zinb(V)$, is unique up to
isomorphisms and can be identified to  $(T(V), \prec)$ where $\prec$ is the half-shuffle
defined as:
\begin{equation*}
v_1 \cdots v_p \prec v_{p+1} \cdots v_n := v_1(v_2 \cdots v_p \Ss v_{p+1} \cdots v_n)
\end{equation*}
\end{proposition}

\subsection{Recall on associative coalgebra}
This section is mainly to fix notations.
\begin{definition}A \emph{coassociative coalgebra}, is a vector space endowed with a
cooperation
$\Delta$ coassociative and counitary : $C \otimes C \stackrel{\Delta}{\to} C$ which
verifies the two following commutative diagrams:
\begin{displaymath}
\xymatrix{
C\ar[r]^{\Delta} \ar[d]_(0.4){\Delta}&C\otimes C\ar[d]^{\Id\otimes \Delta}\\
C\otimes C \ar[r]^(0.4){\Delta\otimes \Id} & C\otimes C\otimes C \ .\\
}
\end{displaymath}
\begin{displaymath}
\xymatrix{
& C \ar[dl]_{\cong} \ar[d]_{\Delta} \ar[dr]^{\cong} & \\
C \otimes \mathbb{K}& C \otimes C \ar[l]_{id \otimes c}\ar[r]^{c \otimes id}& \mathbb{K}
\otimes C
 \ .}
\end{displaymath}
\end{definition}

\begin{definition}
A coalgebra is said to be \emph{connected} if it  verifies the following property:
\begin{equation*}
\begin{array}{l}
H= \bigcup_{r \geq 0}F_{r}H\\
\textrm{where } F_{0}:=\mathbb{K} 1\\
\textrm{and, by induction } F_{r}:=\left\{x \in H \ \vert \  \bar \Delta(x) \in F_{r-1}
\otimes F_{r-1}\right\} \ ,
\end{array}
\end{equation*}
where, $\bar\Delta (x)=\Delta(x)-1\otimes x-x\otimes 1$ .\\
Note that the connectedness only depends on the cooperation and the unit.
\end{definition}

\begin{definition}
A connected coalgebra $C_{0}$ is said to be \emph{free over the vector space $V$} if
there exists a map
$p :C_{0} \to V$ which satisfies the following universal property:

any map $\Phi : C\to V$, where $C$ is a coaugmented connected coalgebra, such that
$\Phi(1)=0$,
extends uniquely in a coalgebra morphism $\tilde \Phi : C\to C_{0}$.
This can be sumed up in the following commutative diagram:
\begin{displaymath}
\xymatrix{C \ar[dr]^{\Phi} \ar[d]_{\tilde \Phi}&\\
C_{0} \ar[r]^{p} & V \ .}
\end{displaymath}
\end{definition}

\begin{definition}
The tensor module $T(V)$ over the vector space $V$ can be endowed with a structure of
coalgebra with the cooperation $\Delta$ defined as:
$$\Delta(v_1\cdots v_n)=\sum_{p=1} ^{n-1} v_1\cdots v_p \otimes v_{p+1}\cdots
v_n+1\otimes v_1\cdots v_n +v_1\cdots v_n\otimes 1$$ with,
$$\begin{array}{rcl}
\Delta(1)&=&1\otimes 1\\
\Delta(v)&=&v\otimes 1+1\otimes v \quad v \in V\\
\end{array}$$
and the counit $c:T(V) \longrightarrow \mathbb K$ is the projection on the first factor
$V$.
\end{definition}

\begin{proposition}
The tensor coalgebra is the free connected coalgebra up to isomorphisms.  $\Box$
\end{proposition}

\subsection{The Associative-Zinbiel bialgebra}
\begin{definition}
An associative-Zinbiel bialgebra is a vector space $\HH$ endowed with a structure of
associative coalgebra $\Delta:\HH
\to\HH\t\HH$, a structure of Zinbiel algebra $\prec:\HH\t\HH\to\HH$, such that the
following compatibility relation is verified:
\begin{equation*}
\Delta(x \prec y)=x_1\prec y_1 \otimes x_2 \ast y_2, \ \forall  x,y \in H ,
\end{equation*}
with the following convention $(1\prec 1)\otimes (x\ast y)=1 \otimes (x\prec y)$.
\end{definition}
\begin{example}
The tensor module $T(V)$ endowed with the deconcatenation and the half-shuffle product is
an associative-Zinbiel bialgebra.
\end{example}

\begin{theorem}\label{thm:As-Zinb}
Let $H$ be an As-Zinb bialgebra over a field $\mathbb{K}$ of any characteristic. The
following are equivalent:
\begin{enumerate}
\item $H$ is connected,
\item $H$ is isomorphic to $(Zinb(V), \prec, \Delta)$ as a bialgebra.
\end{enumerate}
\end{theorem}
This theorem can now be seen as a particular case of the structure theorem for associative-dendriform bialgebras done by M. Ronco in \cite{R}. To do so one must rephrase her article in terms of generalised bialgebra theory and realise that a Zinbiel algebra is a kind of commutative Dendriform algebra. Then, one can show that the primitive structure found in the Associative-Dendriform case crushes to a vector space structure, \cite{B}. 
We give in this paper a straightforward proof of the theorem.
\subsection{Proof of the theorem}
\begin{definition}
The \emph{convolution} of two Zinbiel algebra morphisms  $f$ and  $g$ are defined by:
\begin{equation*}
f \star g := \prec \circ (f \otimes g) \circ \Delta  \ .
\end{equation*}
Note that this convolution is not associative.
\end{definition}
\begin{lemma}
Let $H:=\KK\oplus\bar H $ be an Zinbiel-associative bialgebra. The map $e: H
\longrightarrow H$ is defined:
\begin{equation}
e:= J-J \star J+ (J\star J) \star J - \cdots +(-1)^{n-1} \ \GS \star n J+ \cdots
\end{equation}
where $J:=Id-uc$ and $\GS \star n J:= (\cdots((J\star J)\star J)\cdots\star J)$,
satisfied the following properties:
\begin{enumerate}
\item  $\Im e = \Prim H $,
\item  $\forall x,y \in \bar H$,  $ \ e(x \prec y)=0$,
\item  $e$ is an idempotent,
\item for $H=(Zinb(V)_+,\prec,\Delta)$ defined above, $e$ is identical on $V$
and zero on the other components.
\end{enumerate}
\end{lemma}
\begin{proof}
Note that $e=\Id-\mu\circ\DD$. The first assertion is done by induction on the degree of $x\in F_{r}\bar H$. The second assertion is proven by the bialgebra compatibility relation. The third assertion is obtained by direct computation taking into account  the second assertion. The last assertion is done by direct computation with the second assertion.
\end{proof}
\begin{definition}
Let $\PBT_n$ define the set of rooted trees with $n$ leaves.
We define the operations in the free magmatic algebra $Mag(V):=\KK[\PBT_n]\t V^{\t n}$,
for all $T \in Mag(V)$, \begin{eqnarray*}
T\G n:=(T\cdot (T\cdots (T\cdot (T\cdot T)))) \\
\G nT:=((((T\cdot T)\cdot T) \cdots )\cdot T)
\end{eqnarray*}
We define the completion of the magmatic algebra  $Mag(\mathbb K)^{\wedge}$ as
$Mag(\mathbb K)^{\wedge}~=~\prod_{n\geq 0}\KK[\PBT_n]$, where we denote the first
generator $|$ by $t$. This allows to define formal series in $Mag(\mathbb K)^{\wedge}$.
\end{definition}
\begin{proposition}\label{prop:tree}
In $Mag(\mathbb K)^{\wedge}$, the formal series \begin{eqnarray*}
g(t)&=& t-\G 2 t+\G 3 t+\cdots +(-1)^{n+1}\ \G n t+\cdots \\
f(t)&=&t+t\G 2+t\G 3+\cdots+t\G n+\cdots \ ,
\end{eqnarray*}
are inverse for the composition.
\end{proposition} \begin{proof}
The proof is done by induction. Direct calculation shows that up to rank 1 the property
is verified.
 Suppose that the property is verified up to rank $n$, then:
\begin{eqnarray*}
(f\circ g(t))_{n+1}&=& \sum_{i_1+ \cdots +i_q=n+1}(-1)^{n-q}(\cdots((t^{i_1}\cdot
t^{i_2}) \cdot t^{i_3})\cdot\cdots \cdot t^{i_q})\\
&=&\sum_{i_q} (-1)^{n-q}\Big(\underbrace{\sum_{i_1+ \cdots +i_{q-1}=n+1-i_q}
(-1)^{n-q}
(\cdots((t^{i_1}\cdot t^{i_2}) \cdot \cdots \cdot t^{i_{q-1}})}_{\textrm{by
induction}=0}\Big)\cdot t^{i_q})\\
\end{eqnarray*}
We verify that the right inverse is a left inverse too, as in the associative context:
\begin{eqnarray*}
f^{-1}=f^{-1}\circ (f\circ g)=(f^{-1}\circ f)\circ g = g \ .
\end{eqnarray*}
Therefore, we proved  $f\circ g=Id$ and $g\circ f=Id$.
\end{proof}
\begin{proof}[Proof of theorem \ref{thm:As-Zinb}] We denote $V:=\Prim H$. We define the
map $$G:\overline{H} \longrightarrow Zinb(V):x\mapsto J(x)-\sum (-1)^{n-1}\GS \star n J \
,$$
where, $\GS \star n J:= (((J\star J)\star\cdots\star J)\star J)$ and the map
$$F:Zinb(V)  \longrightarrow \overline{H}:x \mapsto J(x) + \sum J^{\star n} \ ,$$
where, $J^{\star n}:= (J\star (J\star(\cdots\star (J\star J))))$.
Moreover, we define the two formal series in $Mag(V)^{\wedge}$ :
\begin{eqnarray*}
g(t)&=& t-\G 2 t+\G 3 t+\cdots +(-1)^{n+1}\ \G n t+\cdots \\
f(t)&=&t+t^2+t\cdot t^2+ \cdots + (t\cdot (t\cdot (\cdots (t\cdot t^2))))+ \cdots\\
&=&t+t\G 2+t\G 3+\cdots+t\G n+\cdots \ ,
\end{eqnarray*}
which are inverse for the composition by proposition \ref{prop:tree}.
We apply these series on $Hom_\mathbb{K}(H,H)$ sending $1$ on $0$ using $\star$ as
multiplication,
thanks to the following map:
\begin{equation*}
\begin{array}{rcl}
Mag(\mathbb K)^{\wedge}& \longrightarrow & Hom_\mathbb{K}(H,H)\\
t&\mapsto & J\\
\phi(t)=\sum a_n t\G n&\mapsto &\phi^\star(J)=\Phi(x)=\sum a_nJ^{\star n}(x)\\
\psi(t)=\sum b_n \G n t&\mapsto &\psi^\star(J)=\Psi(x)=\sum b_n \GS \star n J(x)\\
\phi\circ\psi(t)&\mapsto & (\phi\circ\psi)^\star(J)
=\Phi\circ\Psi(x)=\phi^\star(J)\circ\psi^\star(J) \end{array}
\end{equation*}
We obtain $e=g^\star J $ and \begin{eqnarray*}
F \circ G =f^{\star} \circ g^{\star}(J)= (f \circ g)^{\star}(J)=Id^{\star}(J)=J\ ,\\ G
\circ F =g^{\star} \circ f^{\star}(J)= (g \circ f)^{\star}(J)=Id^{\star}(J)=J\ .
\end{eqnarray*}
This ends the proof as $J=Id$ on $\bar H$. \end{proof}

\section{The Zinbiel-Associative structure theorem}\label{App:Zinb-As}
This section is just a dualisation of the above section.
\begin{definition}
A \emph{Zinbiel} coalgebra is a vector space $C$ endowed with a co-operation
$\DD:C\rightarrow C\t C$ such that :
$$
(\DD\t\Id)\circ\DD=(\Id\t\DD)\circ\DD+(\Id\t\tau\DD)\circ\DD \ ,
$$
where $\tau:C\t C\to C\t C$ is the map which interchanges the two factors: $\tau(x\t
y)=y\t x$.

A Zinbiel coalgebra is said to be \emph{counital} if it admits a linear map
$c:C\rightarrow \KK$ such that :
$$
\Bigg\{ \begin{array}{lcl}
(c\t\Id)\circ\DD &=& 0 \ ,\\
(\Id\t c)\circ\DD&=& \Id \ . \end{array} $$
It is to be noted that $(c\t c)\circ\DD$ is not defined.
This notion is dual to the notion of Zinbiel algebra (originally called dual Leibniz
algebra in \cite{LZin}) \end{definition}
\begin{remark}The co-operation $\Delta :=\tau\DD+\DD : C \rightarrow C\t C $
is coassociative co-commutative and counital.
\end{remark}
\begin{definition}
A connected coalgebra $H=\KK\oplus \bar H$ is a coalgebra verifying the following
property:
\begin{equation*}
\begin{array}{l}
H= \bigcup_{r \geq 0}F_{r}H\\
\textrm{where } F_{0}:=\mathbb{K} 1\\
\textrm{and by induction } F_{r}:=\left\{x \in H \ \vert \  \bar\DD(x) \in F_{r-1}
\otimes F_{r-1}\right\} \ .
\end{array}
\end{equation*}
where $\bar\DD (x)=\DD(x)-x\otimes 1$ .
\end{definition}

\begin{example}
A \emph{co-shuffle} coproduct can be defined on the tensor module 
$T(V)$ over a vector space $V$ as follows :\begin{equation*}
\SS(v_1\cdots v_p v_{p+1} \cdots v_n) := \sum_{p+q=n}\sum_{\underline i \in
Sh_{p,q}}v_{i_1}\ldots v_{i_p}\t v_{i_{p+1}}\ldots v_{i_n} \in V^{\otimes n} \ ,
\end{equation*}
where the sum is extended over all $(p,q)$-shuffles $\underline i$ (i.e. in the
multi-index $\underline i=(i_1,\ldots,i_n)$ the integers $1, \ldots, p$ are ordered and
so are $p+1, \ldots, n$).

The tensor module $T(V)$ endowed with the \emph{co-half shuffle} $\DD:=\Id\t\SS$ is the
cofree Zinbiel coalgebra.
\end{example}

\begin{definition}
A \emph{Zinbiel-associative bialgebra} $H=(H,\mu, \DD)$ is a vector space $H=\bar H
\oplus  \mathbb{K}\ 1$ endowed with a co-unital Zinbiel co-operation $\DD$ and an
associative operation $\mu$ verifying the following compatibility relation :
$$
\DD\circ\mu=(\mu\t\mu)\circ(\Id\t\tau\t\Id)\circ(\DD\t\Delta) \ ,
$$
\end{definition}

\begin{example}
The tensor module endowed with the concatenation product $\cdot$ and the co-half shuffle
$\DD$ is a Zinbiel-associative bialgebra. \end{example}

\begin{theorem}{\cite{B}}\label{theo:Zinb-as}
Let $\HH$ be a Zinbiel-associative bialgebra over the field $\mathbb{K}$ without any
assumption on its characteristic. The following are equivalent :
\begin{enumerate}
\item $H$ is connected,
\item $H$ is isomorphic to $(T(V), \cdot, \DD)$.
\end{enumerate}
\end{theorem}

\section{Leibniz homology and the Zinbiel coalgebra structure \cite{Lleib}}\label{App:Leib-Zinb}
Let $V$ be a vector space.
Let $\DD$ denote the co-half shuffle defined on the tensor module $T(V)$.
\begin{proposition}\label{prop:Interchange_d_D}
Let $\DD^{p,q}$ denote the projection of $\DD$ on the vector space $\Vt p\t\Vt q$. Let
$d^n$ be the Leibniz differential on $\Vt n$.
Then the following holds :
$$\DD^{p,q}\circ d^{p+q+1}=(d^{p+1}\t\Id)\circ\DD^{p+1,q}+(\Id\t d^{q+1})\circ\DD^{p,q+1}
\ .$$
\end{proposition} The proof is done by dualizing the proof of J.-L. Loday in \cite{Lleib}.
\begin{proof}
First, we shall compute the number of terms appearing on each side of the equation.
On the right hand side there are exactly :
$$\binom{p+q-1}{p-1}\frac{(p+q)(p+q+1)}{2}=\frac{(p+q+1)!}{2(p-1)!q!} $$ terms.
On the left hand side there are :
$$\frac{p(p+1)}{2}\binom{p+q}{q}+\frac{q(q+1)}{2}\binom{p+q}{p-1}=
\frac{(p+q+1)!}{2(p-1)!q!}$$terms.
The number of terms  appearing in each parts of the equation coincide. It suffices
therefore to check that any term on the left side belongs to the set of elements
appearing in the right hand side.

To ease the proof we introduce the following operator $\delta_i^j:\Vt n\rightarrow \Vt
{n-1}$ for $1\leq i<j\leq n$ is defined by : $$\delta_i^j(x_1\ldots
x_n):=x_1\t\ldots\t[x_i,x_j]\t\ldots\t x_n \ ,$$ so that $d^n=\sum_{1\leq i<j\leq
n}(-1)^j\delta_i^j$.

There are two cases to be considered. Let $\sigma^\star$ be a $(p,q)$ co-shuffle.
Consider the element $(\Id\t\delta_k^l)\circ(\Id\t\sigma^\star)$ where $1\leq k<l\leq p$.
This operator is part of $\DD^{p,q}\circ d^{p+q+1}$. Indeed, $(\Id \t\delta_k^l)\circ
(\Id\t\sigma^\star)=(\Id\t\omega^\star)\circ\delta_{\sigma^\star(k)}^{\sigma^\star(l) }$
for a certain $(p-1,q)$-coshuffle $\omega$.

The other case is treated analogously, and this ends the proof. \end{proof}

\begin{corollary}\label{cor:homandprim}
Let $(\HH,d)$ be a graded connected differential Zinbiel-associative bialgebra. Then,

$$H_*(\HH,d)\cong T(H_*(\Prim \HH,d)) \ .$$

The primitive part of the homology of a Zinbiel-associative bialgebra is the tensor
module of the homology of the primitive part of the Zinbiel-associative bialgebra.
\end{corollary}
\begin{proof}
By the above theorem \ref{theo:Zinb-as}, we can restrict ourselves to prove the following
:
$$H_*(T(\Prim \HH,d))=T(H_*(\Prim \HH,d)) \ .$$
And it is well known that two functors $T$ and $H_*$ commute, see \cite{Qa} appendix B.
\end{proof}

\end{document}